\documentclass[12pt]{amsart}
\usepackage[left= 1.1in,top= 1.3in,right= 1.1in,bottom= 1.1in]{geometry}

\usepackage{amsmath}
\usepackage{amsfonts}
\usepackage{amssymb}
\usepackage{amsthm}
\usepackage{latexsym}
\usepackage[all]{xy}
\usepackage{multirow}
\usepackage{graphicx}
\usepackage{epsfig}
\DeclareSymbolFont{AMSb}{U}{msb}{m}{n}
\DeclareMathSymbol{\C}{\mathbin}{AMSb}{"43}

\newtheorem{thm}{Theorem}[section]
\newtheorem{cor}[thm]{Corollary}
\newtheorem{lem}[thm]{Lemma}

\newtheorem{prop}[thm]{Proposition}

\newtheorem{definition}[thm]{Definition}
\newtheorem{ex}[thm]{Example}

\theoremstyle{remark}

\numberwithin{equation}{section}

\newcommand{\Z}{\ensuremath{\mathbb{Z}}}
\newcommand{\Zp}{\ensuremath{\mathbb{Z}_{p}}}
\newcommand{\Q}{\ensuremath{\mathbb{Q}}}
\newcommand{\modp}{\ensuremath{(\text{mod }p)}}

\begin{document}

\title{Surgery description of colored knots}
\author{R.A. Litherland}
\author{Steven D. Wallace}
\address{Louisiana State University}
\email{lither@math.lsu.edu}
\email{wallace@math.lsu.edu}
\date{}


\begin{abstract} 

The pair $(K,\rho)$ consisting of a knot $K\subset S^{3}$ and a surjective map $\rho$ from the knot group onto a dihedral group is said to be a \textit{$p$-colored knot}.  In \cite{Mos}, D. Moskovich conjectures that for any odd prime $p$ there are exactly $p$ equivalence classes of $p$-colored knots up to surgery along unknots in the kernel of the coloring.  We show that there are at most $2p$ equivalence classes.  This is an improvement upon the previous results by Moskovich for $p=3$, and $5$, with no upper bound given in general.  T. Cochran, A. Gerges, and K. Orr, in \cite{CGO}, define invariants of the \textit{surgery equivalence} class of a closed $3$-manifold $M$ in the context of \textit{bordism}.  By taking $M$ to be $0$-framed surgery of $S^{3}$ along $K$ we may define Moskovich's \textit{colored untying invariant} in the same way as the Cochran-Gerges-Orr invariants.  This bordism definition of the colored untying invariant will be then used to establish the upper bound.

\end{abstract}

\maketitle


\section{Introduction}\label{sec:Introduction}
\begin{figure}
	\centering
		\includegraphics[scale = 0.5]{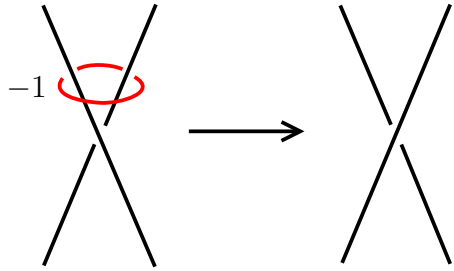}
			\put(-169,65){$-1$}
	\caption{Crossing change due to surgery.}
	\label{fig:CrossingChange}
\end{figure}

It is well known that any knot $K\subset S^{3}$ may be unknotted by a sequence of crossing changes.  A crossing change may be obtained by performing $\pm 1$-framed surgery on $S^{3}$ along an unknot, in the complement of the knot, which loops around both strands of the crossing.  The framing is determined by the sign of the crossing (see Figure \ref{fig:CrossingChange}).  The result of the surgery is once again the $3$-sphere, however the knot $K$ has changed.  By the same token any knot may be obtained from an unknot in the $3$-sphere by $\pm 1$-framed surgery along null-homotopic circles in the complement of the unknot.  This idea is called the \textit{surgery description} of a knot.  For two knots $K_{1}$, and $K_{2}$ in $S^{3}$ we have an equivalence relation defined by $K_{1}\sim K_{2}$ if $K_{2}$ may be obtained from $K_{1}$ via a sequence of $\pm 1$ surgeries along unknots.  Since every knot may be unknotted via this type of surgery we have that $K\sim U$ where $K\in S^{3}$ is any knot and $U\in S^{3}$ is the unknot.  

A $p$-colored knot $(K,\rho)$ is a knot together with a surjection $\rho:\pi_{1}(S^{3}-K)\rightarrow D_{2p}$ from the knot group onto the dihedral group of order $2p$.  An analog of the surgery description of knots for $p$-colored knots is given by restricting the surgeries to those which preserve the existence of a \textit{coloring} $\rho$.  It is natural, then, to ask what the equivalence classes are of $p$-colored knots modulo this surgery relation.  This relation will be refered to as \textit{surgery equivalence in the kernel of $\rho$}, or \textit{surgery equivalence of $p$-colored knots}.

In \cite{Mos}, D. Moskovich proves that for $p=3,5$ there are exactly $p$ equivalence classes.  Moskovich conjectures that this holds for all $p$ and although he has shown that $p$ is a lower bound on the number of equivalence classes in general, no upper bound is given.  In this paper, we will show that the number of surgery equivalence classes of $p$-colored knots is at most $2p$.  More precisely, we will prove the following theorem.

\begin{thm}\label{thm:MainTheorem}\textit{\textbf{Main Theorem}}\\
There are at most $2p$ surgery equivalence classes of $p$-colored knots.  Moreover, if $K_{p}$ denotes the left-handed $(p,2)$-torus knot and $\rho$ is any non-trivial coloring for $K_{p}$ then $$(K_{p},\rho),(K_{p},\rho)\# (K_{p},\rho),\ldots, \#^{p}_{i=1}(K_{p},\rho)$$ are $p$ distinct surgery classes.
\end{thm}

\noindent Note that the list of distinct classes is given in \cite{Mos} but we will use a new definition for his ``colored untying invariant'', denoted $cu(K,\rho)$, to obtain the same result.

\begin{figure}[ht!]
  \vspace{9pt}

  \centerline{\hbox{ \hspace{0.0in} 
    \epsfxsize=2.6in
    \epsffile{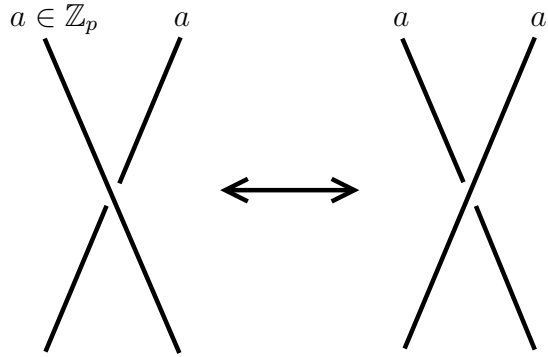}
    	\put(-200,125){$a\in \Z_{p}$}
    	\put(-3,125){$a$}
    	\put(-55,125){$a$}
    	\put(-138,125){$a$}   	
    }
  }

  \vspace{9pt}
  \hbox{\hspace{2.3in} (a) The RR move.} 
  \vspace{9pt}

  \centerline{\hbox{ \hspace{0.0in}
    \epsfxsize=2.6in
    \epsffile{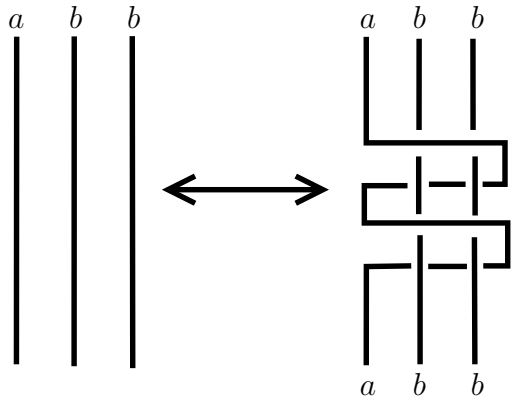}
    	\put(-190,130){$a$}
    	\put(-145,130){$b$}
    	\put(-167,130){$b$}
    	\put(-57,130){$a$}
    	\put(-37,130){$b$}
    	\put(-15,130){$b$}
    	\put(-57,-9){$a$}
    	\put(-37,-9){$b$}
    	\put(-15,-9){$b$} 
    }
  }

  \vspace{9pt}
  \hbox{\hspace{2.3in} (b) The R2G move.} 
  \vspace{9pt}

  \caption{The RR and R2G moves.}
  \label{fig:RRmove-R2Gmove}
\end{figure}

One way to attempt to establish an upper bound on the number of surgery equivalence classes is by using some basic moves on diagrams which preserve colorability and thereby perhaps reducing the crossing number of the diagram or knot.  This is a direct analog to the classical unknotting result where the basic move is a simple crossing change.  It was in this way that Moskovich proved his result for $p=3,5$.  These basic moves are called the $RR$ and $R2G$-moves shown in Figure \ref{fig:RRmove-R2Gmove} (a) and (b).
 
Another interesting question arises: Is there always a finite list of basic moves which are sufficient to describe surgery equivalence of colored knots as Reidemeister moves do for isotopy of knots?  Although it is not proven directly in his paper, Moskovich's result for $p=3,5$ gives a sufficient list of moves which may be used to untie a colored knot consisting of the $RR$, and $R2G$-moves, along with the ``unlinking of bands.''  So the answer is yes for $p=3,5$ but it is unknown otherwise.  The following example shows a non-trivial relation between $3$-colored knots.

\begin{figure}[hb]
	\centering
		\includegraphics[scale = 0.7]{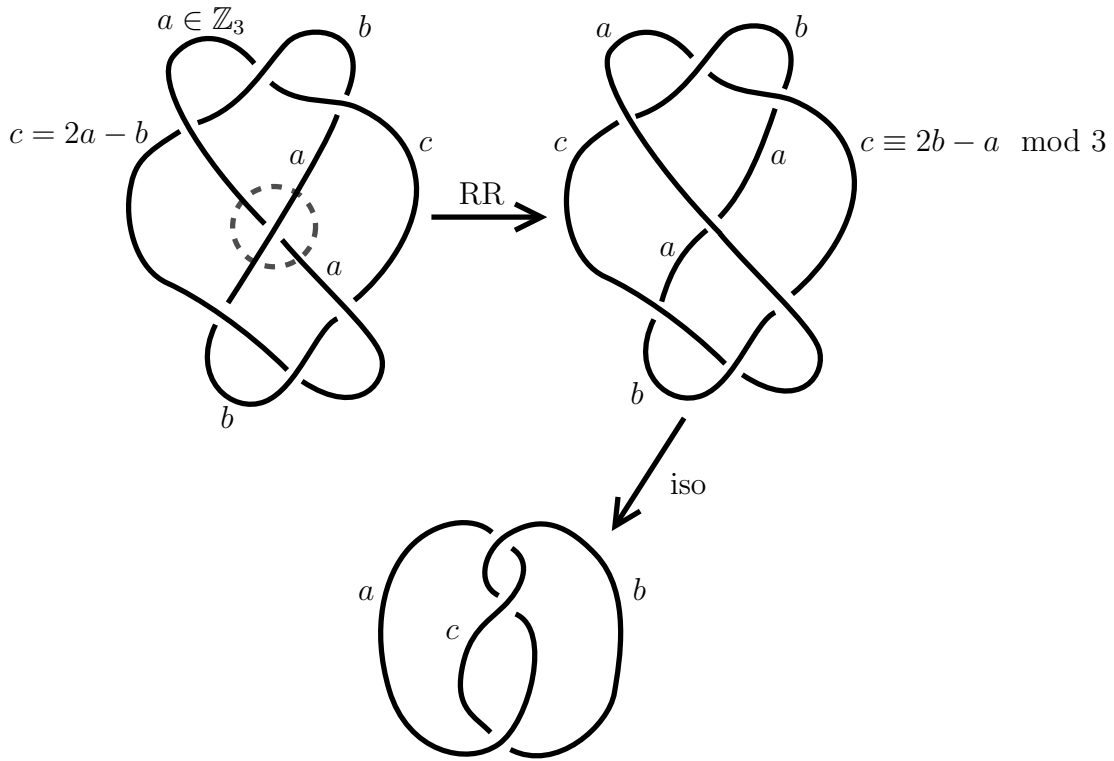}
			\put(-284,276){$a\in \Z_{3}$}
			\put(-234,225){$a$}
			\put(-220,183){$a$}
			\put(-208,273){$b$}
			\put(-260,125){$b$}
			\put(-340,233){$c=2a-b$}
			\put(-185,230){$c$}
			\put(-118,273){$a$}
			\put(-52,225){$a$}
			\put(-94,190){$a$}
			\put(-43,273){$b$}
			\put(-105,134){$b$}
			\put(-18,230){$c\equiv 2b-a \mod{3}$}
			\put(-134,230){$c$}
			\put(-208,60){$a$}
			\put(-104,60){$b$}
			\put(-175,45){$c$}
			\put(-170,210){RR}
			\put(-90,100){iso}
			
	\caption{The $7_{4}$ knot is surgery equivalent to the trefoil knot.}
	\label{fig:7_4vsTrefoil}
\end{figure}

\begin{ex}\label{ex:3_1vs7_4}
The right-handed trefoil knot ($3_{1}$) and the $7_{4}$-knot are surgery equivalent $p$-colored knots.
\end{ex}
\begin{proof}
Performing a single $RR$-move changes $7_{4}$ into the trefoil as in Figure \ref{fig:7_4vsTrefoil}.  Note that this also shows that the mirror images of these knots are equivalent.  However, neither of these knots is surgery equivalent to its mirror image.  This may be seen by calculation of the \textit{colored untying invariant} as in Section \ref{sec:cuEx}.  
\end{proof}

We will not attempt a direct proof of Theorem \ref{thm:MainTheorem} for $p>5$ as Moskovich does for the first two cases.  Instead we will show that an analog to the \textit{Lickorish-Wallace Theorem} and some basic bordism theory suffices to show that there are no more than $2p$ classes.  The Lickorish-Wallace Theorem states that any closed, oriented, connected $3$-manifold may be obtained by performing Dehn surgery on a link in $S^{3}$ with $\pm 1$-framings on each component.  Furthermore, each component may be assumed to be unknotted.  In \cite{CGO}, T. Cochran, A. Gerges, and K. Orr ask what the equivalence classes of $3$-manifolds would be if we restrict the surgeries to a smaller class of links.  Surgery equivalence of $p$-colored knots may be described in a similar way.

The proof of the Main Theorem then is outlined in four steps.  Step $1$ is to establish a $3$-manifold bordism invariant which coincides with colored knot surgery.  Step $2$ is to show that if two colored knots have bordant knot exteriors with the property that the boundary of the bordism $4$-manifold $W$ is $$\partial W = (S^{3}-K_{1})\coprod (S^{3}-K_{2})\cup_{T^{2}\coprod T^{2}}(T^{2}\times [0,1]),$$ where $(K_{i},\rho_{i})$ are the colored knots, then the colored knots are surgery equivalent.  That is, under these conditions, the bordism may be obtained by adding $2$-handles to the $4$-manifold $(S^{3}-K_{1})\times [0,1]$.  We may then ``fill in'' the boundary of the bordism $4$-manifold by gluing in a solid torus crossed with $[0,1]$.  This new $4$-manifold is a bordism between two copies of $S^{3}$ which corresponds to some surgery description for the $3$-sphere.  So step $3$ is to apply \textit{Kirby's Theorem} to unknot and unlink the surgery curves which may be done by only \textit{handle slides} and \textit{blow-ups} (see \cite{GomSt}).  This establishes a surgery equivalence for the knots that are ``taken along for the ride'' during the handle slides and are unchanged (up to surgery equivalence) by blow-ups.  The final step is to show that if any three colored knots have bordant knot complements, then at least two of the colored knots must be surgery equivalent.

The paper is organized as follows.  First we will precisely state what is meant by $p$-colored knots and surgery equivalence.  Then we will define some invariants of $p$-colored knot surgery equivalence in Section \ref{sec:Invariants}.  There are three types of invariants: the \textit{colored untying invariant}, and the \textit{closed} and \textit{relative bordism invariants}.  The colored untying invariant may be computed using the Seifert matrix as in \cite{Mos}, but we show in Section \ref{sec:GoeritzDef} that it may be defined using the \textit{Goeritz matrix} which allows for a simple and geometric proof of invariance under surgery.  Then, in Section \ref{sec:cuEx}, we compute some examples using the Goeritz definition of the colored untying invariant and hence establishing the lower bound of $p$ for the number of surgery equivalence classes which was previously done using the Seifert matrix.  In Section \ref{sec:BordismInvariants} we show that every one of the colored knot surgery invariants give the same information and thus they are all computable given a diagram for the colored knot.  The bordism invariant $\omega_{2}(K,\rho)$ from Section \ref{sec:ClosedBordInv} is used to relate the colored untying invariant to $\omega_{0}$ which in turn relates the relative bordism invariant $\omega$.  In Section \ref{sec:SurgEquiv}, we will show that a \textit{relative bordism} over the \textit{Eilenberg-Maclane space pair} $(K(D_{2p},1),K(\Z_{2},1))$ between two colored knot exteriors establishes a surgery equivalence between the colored knots $(K_{i},\rho_{i})$ at least half of the time.  This gives an upper bound on the number of equivalence classes for any $p$ which is the main result of the paper.

\section{Colored knots} \label{sec:ColoredKnots}

We will first introduce what is meant by a $p$-colored knot and surgery equivalence of $p$-colored knots.

\subsection{Definitions}

Throughout, let $p$ denote an odd prime.  

\begin{definition}
The pair $\left(K,\rho\right)$ consisting of a knot $K\subset S^{3}$ and a surjective homomorphism, $\rho : \pi_{1}(S^{3}-K, x_{0})\rightarrow D_{2p}$, from the knot group with basepoint $x_{0}$ onto the dihedral group of order $2p$, up to an inner automorphism of $D_{2p}$, is said to be a \textit{$p$-colored knot}.  The knot $K$ is said to be \textit{$p$-colorable} with \textit{coloring} given by $\rho$.
\end{definition}

A coloring $\rho$ is only considered up to an inner automorphism of the dihedral group.  In particular, this means that two $p$-colored knots $\left(K_{1},\rho_{1}\right)$ and $\left(K_{2}, \rho_{2}\right)$ are in the same \textit{coloring class} if $K_{1}$ is ambient isotopic to $K_{2}$ and that the following diagram commutes:
\begin{equation*}
\xymatrix{
	\pi_{1}(S^{3}-K_{1},x_{1}) \ar[d]^{\epsilon} \ar[r]^<<<{\rho_{1}} & D_{2p}\ar[d]^{\sigma} \\
	\pi_{1}(S^{3}-K_{2},x_{2}) \ar[r]^<<<{\rho_{2}} & D_{2p}    
}
\end{equation*}
\noindent where, $\sigma : D_{2p}\rightarrow D_{2p}$ is an inner automorphism and $\epsilon: \pi_{1}(S^{3}-K_{1}, x_{1}) \rightarrow \pi_{1}(S^{3}-K_{2}, x_{2})$ is the isomorphism given by $$\left[\alpha\right] \in \pi_{1}(S^{3}-K_{1}, x_{1}) \mapsto \left[h^{-1}\alpha h\right]=\left[h\right]^{-1}\left[\alpha\right]\left[h\right] \in \pi_{1}(S^{3}-K_{2}, x_{2})$$ where $h$ is any fixed path from $x_{2}$ to $x_{1}$ in $S^{3}-K_{1}$.

If we let $K_{1}=K_{2}$ we see that the choice of a different basepoint results in an inner automorphism of the knot group and thus results in an inner automorphism of the dihedral group.  So the definition is well-defined for any choice of basepoint.  We will then ignore basepoints from now on and denote a coloring simply by a surjection $$\rho : \pi_{1}(S^{3}-K)\rightarrow D_{2p}$$ from the knot group onto the dihedral group.

\begin{definition}
Two $p$-colored knots $(K_{1},\rho_{1})$ and $(K_{2},\rho_{2})$ are \textit{surgery equivalent in the kernel of $\rho$} (or simply \textit{surgery equivalent}) if $K_{2}\in S^{3}$ may be obtained from $K_{1}$ via a sequence of $\pm 1$-framed surgeries of $S^{3}$ along unknots in the kernel of $\rho_{1}$.  Furthermore, $\rho_{2}$ must be compatible with the result on $\rho_{1}$ after the surgeries.  That is, if $K(D_{2p},1)$ denotes an \textit{Eilenberg-Maclane space} over the dihedral group then
\begin{equation*}
\xymatrix{
	S^{3}-K_{1} \ar[d]_{\nu} \ar[dr]^{f_{1}} & \; \\
	S^{3}-K_{2} \ar[r]^<<<{f_{2}} & K(D_{2p},1)    
}
\end{equation*}
\noindent is a commutative diagram where $\rho_{i}$ are the induced maps of the $f_{i}$ on $\pi_{1}$ and $\nu$ is the map resulting from surgery restricted to $S^{3}-K_{1}$.
\end{definition}

So there are two conditions for surgery equivalence of $p$-colored knots: $(1)$ the knots must be surgery equivalent in the classical sense with the restriction that the surgery curves are in the kernel of the coloring, and $(2)$ the coloring of the second knot arises from the coloring of the first knot via surgery.  Notice that $(1)$ assumes that the surgery curves are unknotted with $\pm 1$-framings.

 


We will now define what we mean by a \textit{based $p$-colored knot}.

To do this first recall that a \textit{Fox coloring} is classically described by a labeling of the arcs in a diagram for $K$ with the ``colors'' $\left\{0, \ldots, p-1 \right\}$ (see \cite[Chapter $IV$, Exercise 6]{Cr-Fox}).  At each crossing, the labeling must satisfy the \textit{coloring condition} which requires that the sum of the labels of the underarcs must equal twice the label of the overarc modulo $p$.  We also require that the coloring be \textit{nontrivial}, that is, we require that more than one color is used.  Then such a labeling defines a surjection $\rho:\pi_{1}(S^{3}-K)\rightarrow D_{2p}=\left\langle s, t | s^{2}=t^{p}=stst=1\right\rangle$ by the rule $\rho([\mu])=ts^{l}$ where $l$ is the label given to the arc corresponding to the meridian $\mu$.  Conversely, in a coloring meridians are necessarily mapped to elements of order two in the dihedral group since all meridians are conjugate in the knot group and the coloring is a surjection.  That is a coloring map determines a labeling of any diagram for the knot.

\begin{figure}[hbtp]
  \vspace{9pt}

  \centerline{\hbox{ \hspace{0.0in} 
    \epsfxsize=2.6in
    \epsffile{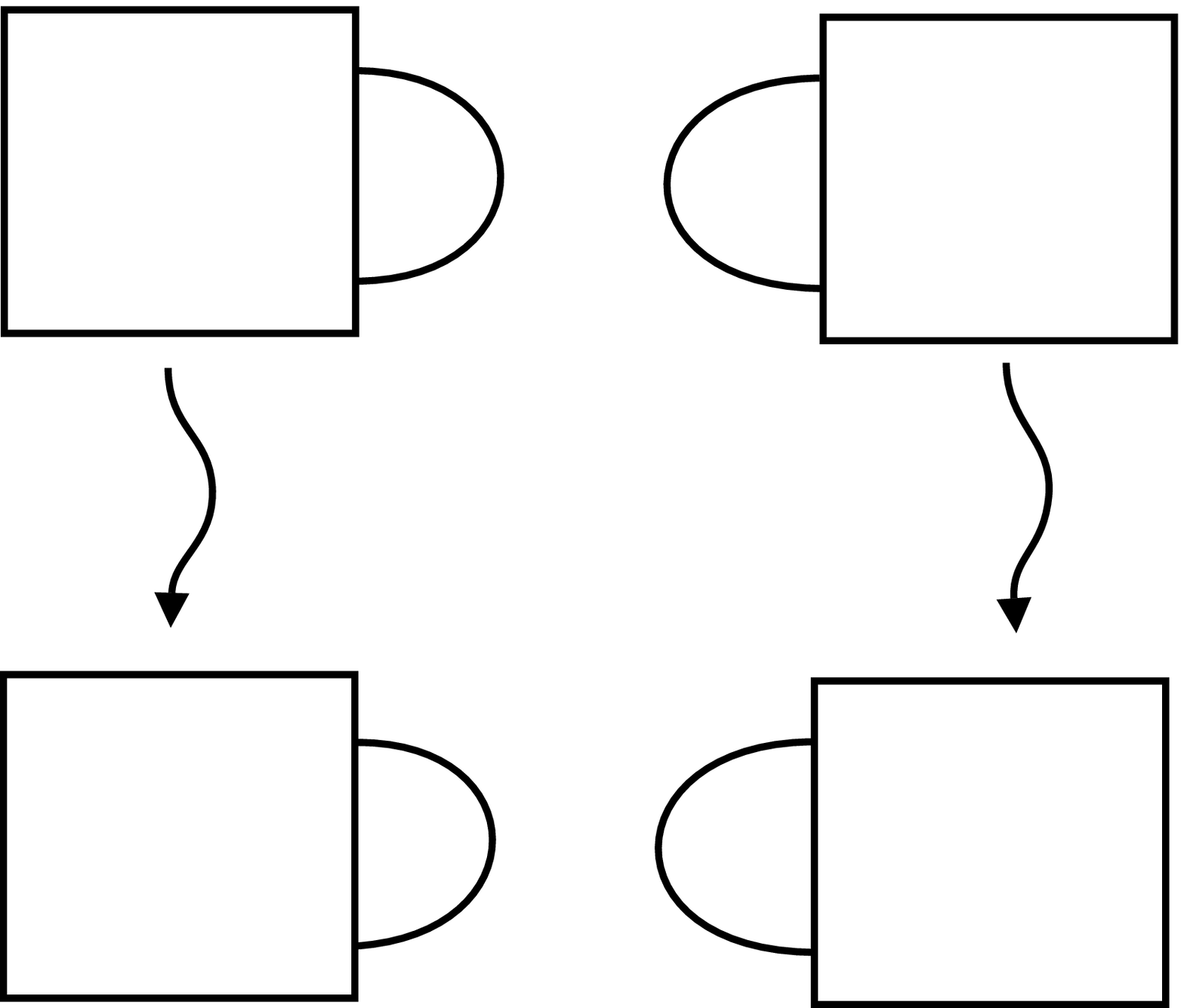}
    	\put(-167,130){$K_{1}$}
    	\put(-38,130){$K_{2}$}
    	\put(-110,145){$a$}
    	\put(-85,145){$a$}
    	\put(-167,24){$K_{1}$}
    	\put(-38,24){$K_{2}$}
    	\put(-85,38){$0$}
    	\put(-110,38){$0$}    	
    }
  }

  \vspace{9pt}
  \hbox{\hspace{1.4in} (a) Relabeling colored knots} 
  \vspace{9pt}

  \centerline{\hbox{ \hspace{0.0in}
    \epsfxsize=2.6in
    \epsffile{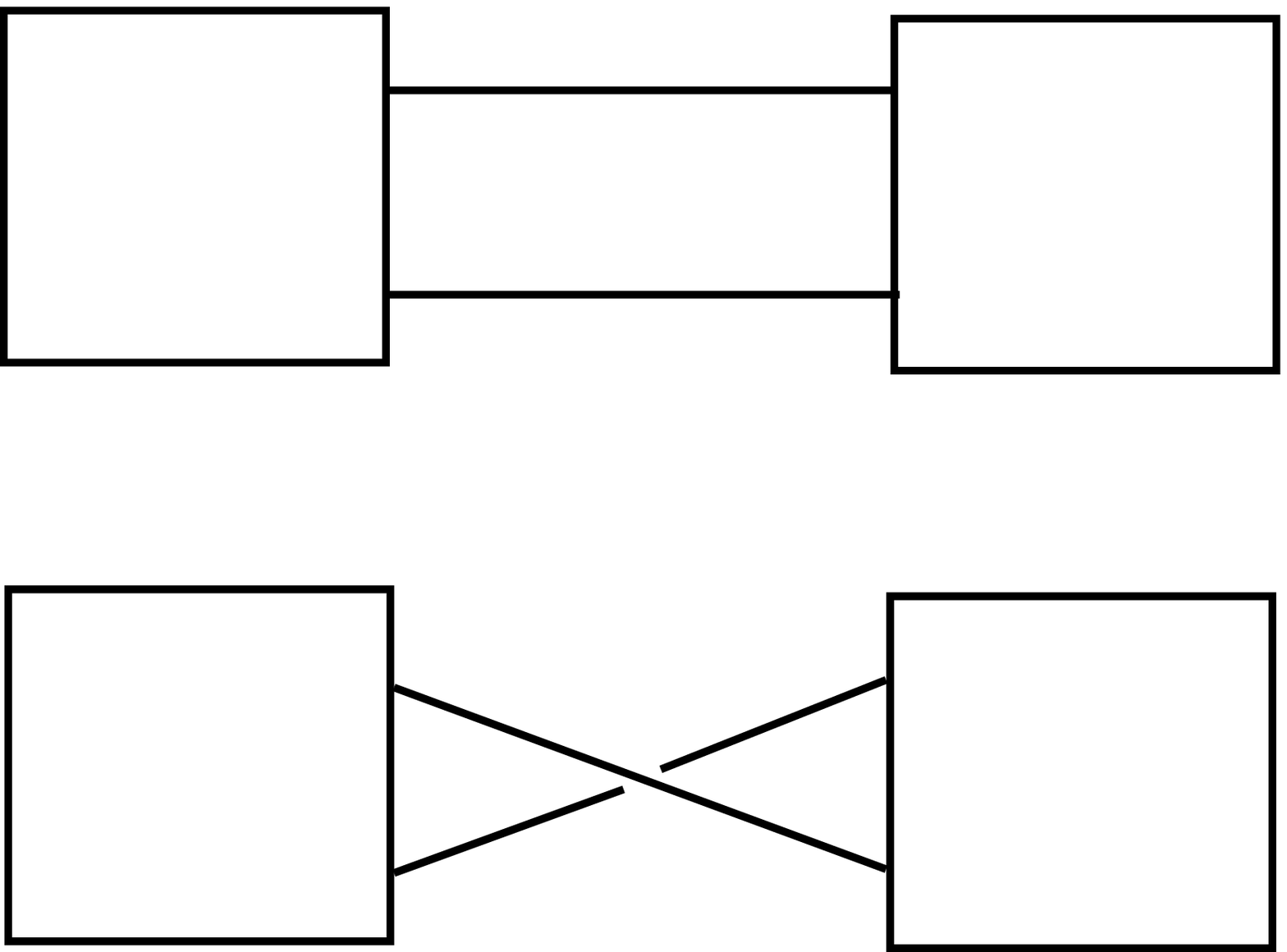}
    	\put(-167,107){$K_{1}$}
    	\put(-38,107){$K_{2}$}
    	\put(-95,130){$0$}
    	\put(-95,100){$0$}
    	\put(-185,67){\large{or (depending on orientation)}}
    	\put(-167,24){$K_{1}$}
    	\put(-38,24){$K_{2}$}
    	\put(-80,36){$0$}
    	\put(-115,36){$0$} 
    }
  }

  \vspace{9pt}
  \hbox{\hspace{1.4in} (b) Connected sum of colored knots after relabeling} 
  \vspace{9pt}

  \caption{Connected sum of colored knots}
  \label{fig:coloredknot1-2a}
\end{figure}

Since we may alter any coloring by an inner automorphism of $D_{2p}$ we may assume that any one arc we choose in a diagram for $K$ be labeled with the color $0$.  We may assume, then, that for any meridian $m$ of $K$ there is an equivalent coloring $\rho$ which maps $m$ to $ts^{0}\in D_{2p}$.  We call the triple $(K,\rho,m)$ a \textit{based $p$-colored knot}.  Therefore, given $\left(K_{1}, \rho_{1}\right)$ and $\left(K_{2},\rho_{2}\right)$ where the $\rho_{i}$ are defined by a nontrivial labeling of a diagram for (oriented) knots $K_{i}$ we may take $(K_{1}$ $\#$ $K_{2},\rho_{3})$ to be the usual connected sum of oriented knots with $\rho_{3}=\tilde{\rho}_{1}\; \# \;\tilde{\rho}_{2}$ (see Figure \ref{fig:coloredknot1-2a}).  Part (a) of the Figure illustrates that we may assume that the $\left(K_{i}, \rho_{i}\right)$ are actually the based $p$-colored knots $\left(K_{i}, \rho_{i},m_{i}\right)$ where $m_{i}$ is the meridian that corresponds to the chosen arc of the diagram for $K_{i}$.

\begin{figure}[hbt]
	\centering
		\includegraphics[scale = 0.5]{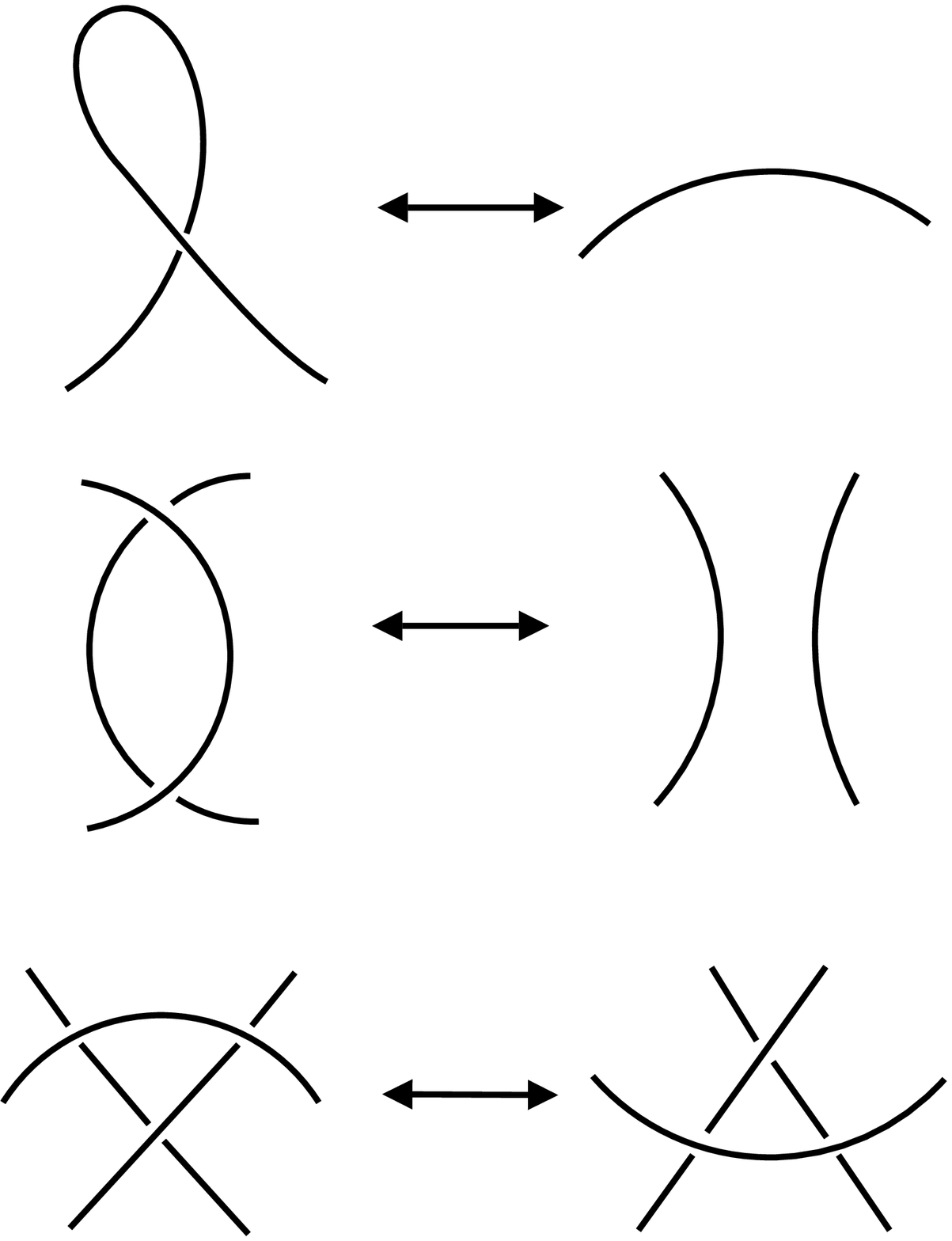}
			\put(-200,300){$a$}
			\put(-80,255){$a$}
			\put(-215,185){$a$}
			\put(-161,185){$b$}
			\put(-240,150){$2a-b$}
			\put(-72,175){$a$}
			\put(-20,175){$b$}
			\put(-225,63){$a$}
			\put(-152,63){$b$}
			\put(-185,63){$c$}
			\put(-236,-4){$2c-b$}
			\put(-174,-4){$2c+a-2b$}
			\put(-25,-4){$2c+a-2b$}
			\put(-85,-4){$2c-b$}
			\put(-62,60){$a$}
			\put(-28,60){$b$}
			\put(-5,30){$c$}
			
	\caption{Colorability is independent of choice of diagram}
	\label{fig:colored_knot3}
\end{figure}

To verify that this process is well-defined for any choice of diagram, we must establish the existence and uniqueness of labelings for each Reidemeister move.  This is done in Figure \ref{fig:colored_knot3}.

\begin{figure}
	\centering
	\includegraphics[scale = 0.60]{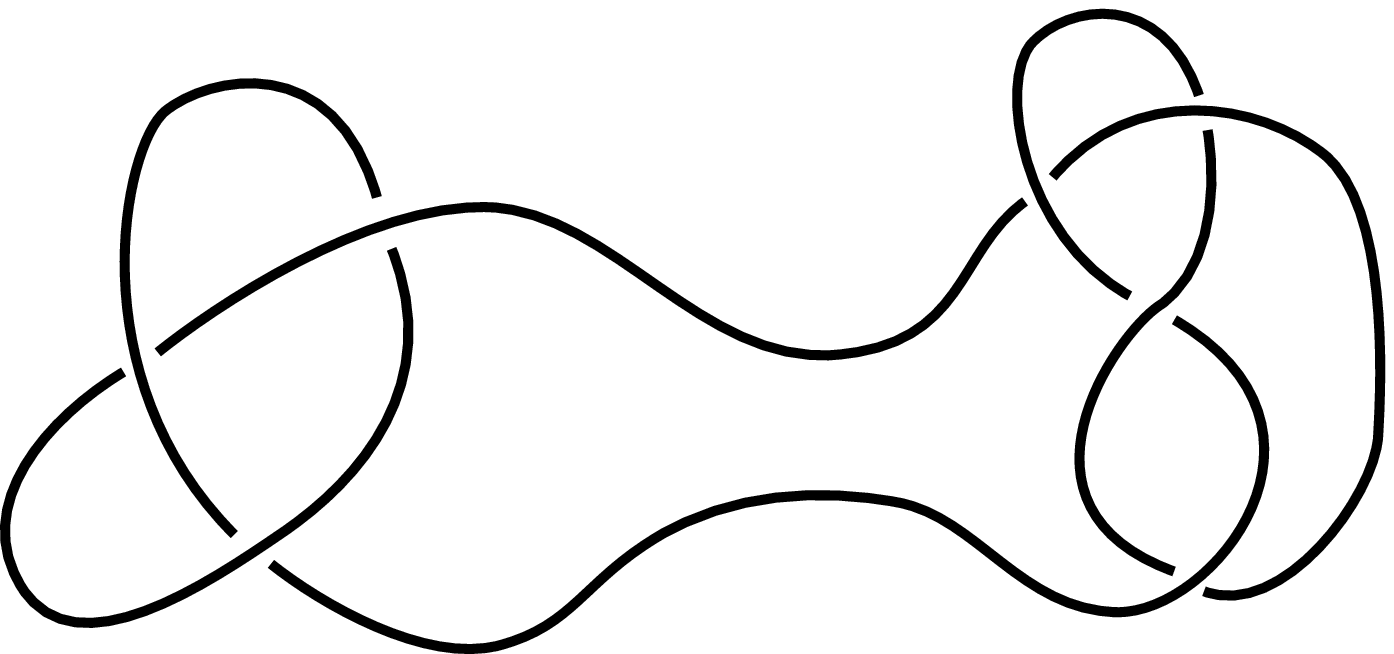}
	\put(-72,92){$0$}
	\put(-100,55){$0$}
	\put(-100,17){$0$}
	\put(-28,78){$0$}
	\put(-2,75){$0$}
	\put(-215,13){$1$}
	\put(-200,103){$2$}
	\caption{A prime $3$-colored knot}
	\label{fig:colored_knot4}
\end{figure}

Unfortunately, the notion of prime $p$-colored knots is slightly different from the usual notion of a prime knot.  For example let $K=K_{1}\; \# \; K_{2}$ where $K_{1}$ is the left-handed trefoil and $K_{2}$ is the figure eight knot.  Then $K$ is $3$-colorable since we can label $K$ using all $3$ colors as in Figure \ref{fig:colored_knot4}.  A knot is $p$-colorable if and only if its determinant is divisible by $p$ \cite{Liv}.  So as $det(K_{2})=5$ and is thus not divisible by $3$, we have that no non-trivial coloring of $K_{2}$ exists.  Therfore, $(K,\rho)\neq (K_{1},\rho_{1})\; \# \; (K_{2}, \rho_{2})$ for any $3$-colorings $\rho_{1}$ and $\rho_{2}$.

\section{Surgery equivalence invariants} \label{sec:Invariants}

As we have seen in Example \ref{ex:3_1vs7_4}, it possible to show that two $p$-colored knots are surgery equivalent directly in some cases.  However, much like trying to distinguish knots by using Reidemeister moves, it is impossible to prove that two $p$-colored knots are not surgery equivalent by simply using a collection of moves on diagrams.  In fact, it is often difficult to show that two knots are the same using Reidemeister moves, and surgery equivalence of $p$-colored knots faces the same type of difficulty.  It is useful then to define algebraic invariants to help distinguish between knot types and the same is true for surgery equivalence.

We do not have a complete list of moves to determine surgery equivalence of $p$-colored knots so we may not simply check an analog to the Reidemeister moves.  Instead, first we must show that the value is unchanged under the choice of $p$-colored knot representative and then we must show that it is invariant under $\pm 1$-surgery.  Then we will show that the three types of $p$-colored knot invariants are in fact three different ways to define the same thing.

\subsection{Preliminaries}

In this section we will introduce some of the background that will be needed in defining the three types of invariants for $p$-colored knots.  The colored untying invariant, $cu$ \cite{Mos}, will arise from the cup product of a certain element $a \in H^{1}(M;\Zp)$ (depending only on the \textit{coloring class}) with its image under the Bockstein homomorphism, $\beta^{1}(a)$ (see \cite{Mun} for a discussion on Bockstein homomorphisms).  In this way we obtain a $\Zp$-valued invariant.  Now we will give a brief overview of the bordism theory needed to define the \textit{closed} bordism invariants $\omega_{2}$, $\omega_{0}$, as well as the \textit{relative} bordism invariant $\omega$.  It will also be useful to recall the definition of the \textit{Goeritz matrix} using the \textit{Gordon-Litherland form} \cite{Gor-Li}.



\begin{figure}[htb]
	\centering
		\includegraphics[scale = .4]{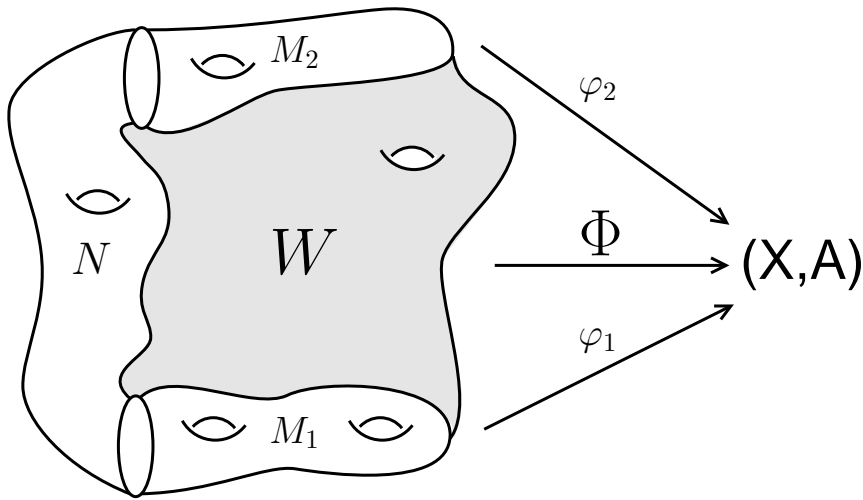}
	\put(-310,85){\Large$N$}
	\put(-235,85){\huge$W$}
	\put(-235,165){\large$M_{2}$}
	\put(-235,21){\large$M_{1}$}
	\put(-118,92){\huge$\Phi$}
	\put(-118,58){\large$\varphi_{1}$}
	\put(-118,153){\large$\varphi_{2}$}
	\caption{Relative bordism over $(X,A)$.}
	\label{fig:invariants1a}
\end{figure}

\begin{definition}
Let $(X,A)$ be a pair of topological spaces $A\subseteq X$.  The $n$\textit{-dimensional oriented relative bordism group of the pair}, denoted $\Omega_{n}(X,A)$, is defined to be the set of \textit{bordism classes} of triples $(M,\partial M, \varphi)$ consisting of a compact, oriented $n$-manifold $M$ with boundary $\partial M$ and a continuous map $\varphi:(M,\partial M)\rightarrow (X, A)$.  The triples $(M_{1},\partial M_{1}, \varphi_{1})$ and $(M_{2},\partial M_{2}, \varphi_{2})$ are in the same bordism class if there exists an $n$-manifold $N$ and a triple $(W, \partial W, \Phi)$ consisting of a compact, oriented $(n+1)$-manifold $W$ with boundary $\partial W = (M_{1} \coprod M_{2})\bigcup_{\partial N} N$ and a continuous map $\Phi:W\rightarrow X$ satisfying $\Phi |_{M_{i}}= \varphi_{i}$ and $\Phi(N)\subseteq A$.  We also require that $M_{1}$ and $M_{2}$ are disjoint and $M_{i}\cap N= \partial M_{i}$ for $i=1,2$.  In this case, we say that $(M_{1}, \partial M_{1}, \varphi_{1})$ and $(M_{2}, \partial M_{2}, \varphi_{2})$ are \textit{bordant over} $(X,A)$ denoted $(M_{1}, \partial M_{1}, \varphi_{1})\sim_{(X,A)} (M_{2}, \partial M_{2}, \varphi_{2})$ (see Figure \ref{fig:invariants1a}).
\end{definition}

\noindent A triple $(M, \partial M, \varphi)$ is \textit{null-bordant}, or \textit{bords}, over $(X,A)$ if it bounds $(W, \partial W, \Phi)$.  That is, it bords if it is bordant to the empty set $\emptyset$.  The set $\Omega_{n}(X,A)$ forms a group with the operation of disjoint union and identity element $\emptyset$.  We will denote $\Omega_{n}(X,\emptyset)$ by $\Omega_{n}(X)$ and so our definition makes sense for pairs $(M,\varphi)=(M,\emptyset,\varphi)$ with $M$ a closed $n$-manifold.

We will only be interested in the case when $n=3$.  In this case, the \textit{Atiyah-Hirzebruch spectral sequence} (see \cite{Whi} for the \textit{extraordinary homology theory} made up of the bordism groups $\Omega_{n}(X, A)$) implies that $\Omega_{3}(X,A)\cong H_{3}(X,A;\Omega_{0}) \cong H_{3}(X,A)$ where $\Omega_{0}\cong \mathbb{Z}$ is the $0$-dimensional bordism group of a single point.  The isomorphism is given by $(M,\partial M, \varphi)\mapsto \varphi_{\ast}(\left[M,\partial M \right])$ where $\left[M, \partial M \right]$ is the fundamental class in $H_{3}(M, \partial M)$.  Furthermore, if we assume that $X$ is an Eilenberg-Maclane space $K(G,1)$ and $A$ is the subspace corresponding to a subgroup $H\subset G$, then the bordism group is isomorphic to the homology of the group $G$ relative the subgroup $H$, that is $\Omega_{3}(X,A)\cong H_{3}(G,H)$.  When $K(H,1)$ is a subspace of $K(G,1)$ we will denote $\Omega_{3}(K(G,1),K(H,1))$ by $\Omega_{3}(G,H)$.

Now for a brief discussion on the Goeritz matrix.  Given a spanning surface $F$ for a link $K$ and a basis $x_{i}$ for its homology, the \textit{Goeritz matrix} is given by evaluating the \textit{Gordon-Litherland form}, $\mathfrak{G}_{F}:H_{1}(F)\times H_{1}(F)\rightarrow \Z$, on the basis elements (see \cite{Gor-Li}).  That is, $G=(g_{ij})$ is defined by $$g_{ij}=\mathfrak{G}_{F}(x_{i},x_{j})=lk(x_{i},\tau^{-1}(x_{j}))$$ where $\tau^{-1}(y)$ is $y$ pushed off in ``both directions.''  Precisely, $\tau: \tilde{F}\rightarrow F$ is the \textit{orientable double covering space} of $F$ (see Chapter $7$ of \cite{Lic}).  Note that $\tilde{F}$ is a connected, orientable surface regardless of the orientability of $F$.

\begin{figure}[htb]
	\centering
		\includegraphics[scale = .7]{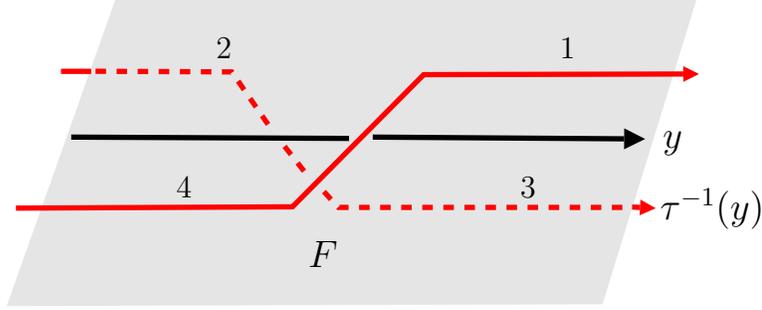}
	\put(-16,60){\large$y$}
	\put(-17,33){\large$\tau^{-1}(y)$}
	\put(-55,94){$1$}
	\put(-185,94){$2$}
	\put(-200,41){$4$}
	\put(-70,41){$3$}
	\put(-150,15){\large$F$}
	\caption{``Double push off'' of a orientation reversing curve $y$.}
	\label{fig:invariants4}
\end{figure}

If $y$ is an orientation preserving loop in $F$ then $\tau^{-1}(y)$ is comprised of two loops, the \textit{positive} $y_{+}$ and \textit{negative} $y_{-}$ push offs on either side of $F$.  Figure \ref{fig:invariants4} illustrates the non-orientable case.  If rather $y$ is orientation reversing then $\tau^{-1}(y)$ is a single loop which double covers $y$.  In this case you can think of $\tau^{-1}(y)$ to be the loop that arises from pushing $y$ off to one side which then comes back around on the other side and vice versa.

\begin{figure}
	\centering
		\includegraphics[scale = .9]{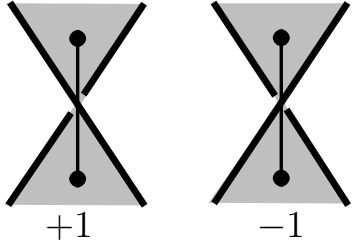}
	\put(-117,-10){\large$+1$}
	\put(-37,-10){\large$-1$}
	\caption{Incidence number at a crossing}
	\label{fig:invariants5}
\end{figure}

Note that the Goeritz matrix for a knot $K$ may also be calculated from a checkerboard coloring for a diagram for the knot (see Chapter $9$ of \cite{Lic}).  First we must pick a white region, the so-called \textit{infinite region} $R_{0}$, and then we number the other white regions $R_{1}, \ldots, R_{n}$.  We then define an incidence number $\iota(c)=\pm 1$ assigned to any crossing $c$ by the rule in Figure \ref{fig:invariants5}.  We define a $(n+1)\times (n+1)$ matrix $(g_{ij})$ for $i \neq j$ by $$g_{ij}=\sum \iota(c),$$ where the sum is over all crossings which are incident with both $R_{i}$ and $R_{j}$.  The diagonal terms are chosen so that the rows and columns sum to $0$, namely $$g_{ii}=-\sum_{l\neq k} g_{lk}.$$  The Goeritz matrix is then obtained from the ``pre-Goeritz matrix'' $(g_{ij})$ by deleting the row and column corresponding to the infinite region.  The group that this matrix presents is independent of the choice of infinite region.

We will use this diagramatic way to calculate the Goeritz matrix in Section \ref{sec:cuEx}.

\subsection{The colored untying invariant}

We will now define precisely the $p$-colored knot invariant $cu$, and later define the bordism invariants $\omega_{2}$, $\omega_{0}$, and $\omega$.

\subsubsection{Moskovich's definition} \label{subsec:Mos}

Throughout, let $(K,\rho)$ be a $p$-colored knot with coloring $\rho: \pi_{1}(S^{3}-K)\rightarrow D_{2p}$ where $D_{2p}=\left\langle s,t \; |\; t^{2}=s^{p}=tsts=1 \right\rangle$ is the dihedral group with $2p$ elements.  Also let $\tilde{X}$ denote the $2$-fold cover of $S^{3}$ branched over $K$.  Let $X_{0}$ denote the manifold obtained from $S^{3}$ by performing $0$-framed surgery along $K$.

Consider the following diagram:
\begin{equation}
\xymatrix{
	 H_{1}(S^{3}-F)\ar[drr]^{\overline{\rho}} & \; & \; \\ 
	 \pi_{1}(S^{3}-F)\ar[u]\ar[rr]^{\rho |_{S^{3}-F}}\ar[d] & \; &  \Zp\ar[d] \\  
	 \pi_{1}(S^{3}-K)\ar[rr]^{\rho} \ar[drr]^{l} & \; & D_{2p}\ar[d] \\ 
	 \; & \; & \mathbb{Z}_{2} \\ 
}  \label{ColoringComDiagram} 
\end{equation}
\noindent with the map $l$ defined by $l(x)=lk(x,K) \; (mod \; 2)$.  Note that the coloring map sends meridians to elements of order $2$, in particular, $\rho(\mu_{i})=ts^{k}$ for some $k\in {0,\ldots, p-1}$ where ${\mu_{i}}$ are \textit{Wirtinger generators} for $\pi_{1}(S^{3}-K)$.  Then the lower triangle of Diagram \ref{ColoringComDiagram} commutes by construction.  Furthermore, we see that if $x$ is a loop in $(S^{3}-F)$ then $lk(x,K)\equiv 0 \; (mod \; 2)$ which is enough to establish the commutativity of the rest of Diagram \ref{ColoringComDiagram}.  Indeed, since $x$ is in the complement of the Seifert surface $F$ we may assume that $lk(x,K)=0$.  Notice that the commutativity of the upper triangle of the diagram is immediate since the image of $\rho |_{S^{3}-F}$ is abelian.

Therefore we have established the existence of a map $f:(S^{3}-F)\rightarrow K(\Zp,1)$ from the complement of the surface to an Eilenberg-Maclane space over $\Zp$.  The map $f$ may be extended to the \textit{unbranched $2$-fold cyclic cover} of $S^{3}-K$ denoted $\tilde{Y}$ which is obtained by gluing two copies of $S^{3}-F$ together along two copies of a bicollar $(F-K)\times (-1,1)$ of the interior of the surface.  Call this ``new'' map $f:\tilde{Y}\rightarrow K(\Zp,1)$.  Now we can form the $2$-fold branched cover $\tilde{X}$ by gluing in a solid torus so that the meridian of the solid torus maps to twice the meridian of the torus boundary of $\tilde{Y}$ (see \cite[Chapters 5 and 10]{Rol}).  Since twice a meridian is mapped trivially by $f$ we may extend this map to the $2$-fold branched cover.  Thus we have a map $f:\tilde{X}\rightarrow K(\Zp,1)$ which will be used to associate $\omega_{2}$ with $cu$ later in Section \ref{sec:cuOmega2}.  

We have shown that
\begin{equation}
\xymatrix{
	 H_{1}(\tilde{X})\ar[ddrr]^{\rho'} & \; & \; \\ 
	 \pi_{1}(\tilde{X})\ar[u]\ar[drr]^{f_{\ast}} & \; & \; \\
	 \pi_{1}(\tilde{Y})\ar[u]\ar[rr]^{f_{\ast}}\ar[d] & \; & \Zp \ar[d] \\  
	 \pi_{1}(S^{3}-K)\ar[rr]^{\rho} \ar[drr]^{l} & \; & D_{2p}\ar[d] \\ 
	 \; & \; & \mathbb{Z}_{2} \\ 
}  \label{ColoringComDiagram2} 
\end{equation}
is a commutative diagram.  So the coloring map $\rho$ restricts in the double covering to a map $$\rho':H_{1}(\tilde{X};\Z)\rightarrow \Zp$$ which corresponds to a cohomology class $$a \in H^{1}(\tilde{X};\Zp)\cong Hom(H_{1}(\tilde{X};\Z),\Zp)$$ by the Universal Coefficient Theorem for Cohomology.  The colored untying invariant is defined to be the cup product of $a$ with its image under the Bockstein homomorphism $\beta^{1}:H^{1}(\tilde{X}; \Zp)\rightarrow H^{2}(\tilde{X}; \Z)$.

\begin{definition}\label{def:cuMos}
Given a $p$-colored knot $(K,\rho)$ the \textit{colored untying invariant} of $(K,\rho)$ is
$$cu(K,\rho):=a\cup \beta^{1}a \in H^{3}(\tilde{X};\Zp)$$
which we may think of as an element of $\Zp \cong H^{3}(\tilde{X};\Zp)$.
\end{definition}

\noindent Note that the isomorphism $\Zp \cong H^{3}(\tilde{X};\Zp)$ is given by evaluation on the fundamental class.  

To show that this is actually an invariant of $p$-colored knots we must assert that it is well-defined for any choice of equivalent coloring.  Invariance of the choice of coloring is clear since $cu$ is defined using homology and cohomology groups which are independent of basepoint and conjugacy class in $\pi_{1}(S^{3}-K)$.  To show that $cu$ is a non-trivial invariant we will introduce a way to compute $cu$ by using the Seifert matrix for a given Seifert surface.  It turns out that there is a way to determine the invariant for any spanning surface (including perhaps a non-orientable surface) by using the \textit{Goeritz} matrix.  We will use this definition to establish non-triviality and invariance under $\pm 1$-framed surgery in the kernel of $\rho$.  Note that Moskovich \cite{Mos} gives an alternate proof of the surgery invariance and does not mention the Goeritz definition. 

Let $F$ be a Seifert surface for $K$ with Seifert matrix $S$ with respect to a basis $x_{1}, \ldots, x_{2k}$ of $H_{1}(F)$.  Let $\xi_{1}, \ldots, \xi_{2k}$ be a basis for $H_{1}(S^{3}-F)$ with orientations so that $lk(x_{i},\xi_{j})=\delta_{ij}$.  The proof of the following lemma can be found in \cite{Mos} and will be omitted here.

\begin{lem}\label{lem:cuCalc1}
Let $v:=(v_{1}, \ldots, v_{2k})^{T}\in \Z^{2k}$ be a column vector such that $$v_{i} \; \modp =\overline{\rho}(\xi_{i})$$ for all $i\in 1, \ldots, 2k$.  Then
$$cu(K,\rho)=2\frac{v^{T}\cdot S \cdot v}{p} \; \modp.$$
The vector $v$ is called a \textit{$p$-coloring vector}.
\end{lem}

If $K=(p,2)$ torus knot, then, for a certain choice of $p$-colorings $\rho_{1}$ and $\rho_{2}$, the lemma may be used to show that $cu(K,\rho_{1})\neq cu(K,\rho_{2})$.  We will show this later in Section \ref{sec:cuEx} using the Goeritz definition of the colored untying invariant defined below.
 
\subsubsection{Goeritz definition}\label{sec:GoeritzDef}

We will now extend Lemma \ref{lem:cuCalc1} to any spanning surface for the knot $K$ including perhaps non-orientable surfaces.  We will use this definition for the colored untying invariant to give a geometric proof that it is a surgery equivalence invariant.

\begin{prop}\label{prop:cuGoer}
The colored untying invariant $cu$ may be calculated using the Goeritz matrix for a diagram for $K$.  That is 
\begin{equation}
cu(K,\rho)=\frac{v^{T}\cdot G \cdot v}{p} \; \modp
\label{eq:cuGoer}
\end{equation} 
where $v$ is any $p$-coloring vector and $G$ is the Goeritz matrix.  \label{lem:cuCalc2}
\end{prop}

Since it is clear from Definition \ref{def:cuMos} that $cu(K,\rho)$ is a well-defined invariant of $p$-colored knots, the fact that this also holds for the Goeritz definition is a corollary to Proposition \ref{prop:cuGoer}.  Thus, we will not give a direct proof of well-definedness under the choices of basis for $H_{1}(S^{3}-F)$ or coloring vector here.  We may also assume that $v^{T}Gv \equiv 0 \; \modp$ as is required for the right side of the above equation to make sense.  We will, however, establish well-definedness under the choice of spanning surface, as this is not clear from Lemma \ref{lem:cuCalc1} for a non-orientable spanning surface.
\begin{lem}
The colored untying invariant is independent of the choice of spanning surface. \label{lem:cuSurf}
\end{lem}
\begin{proof}
Spanning surfaces are related by (i) \textit{$S$-equivalence} in the usual sense (see \cite{BFK}), or (ii) addition or deletion of a single twisted band (see Figure \ref{fig:invariants2}).  Note that operation (ii) may perhaps change the orientability of the resulting surface.  We will now show that the right hand side of equation \ref{eq:cuGoer} is unchanged by all three types of moves.

\begin{figure}
	\centering
		\includegraphics[scale = .4]{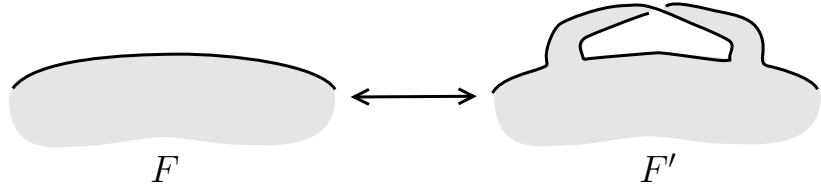}
	\put(-255,-12){\large$F$}
	\put(-70,-12){\large$F'$}
	\caption{Non-orientable $S$-equivalence.}
	\label{fig:invariants2}
\end{figure}

Let $F$ and $F'$ denote $S$-equivalent possibly non-orientable spanning surfaces for $K$ and let $(G,v)$ and $(G',v')$ be the corresponding pairs consisting of a Goeritz matrix and a coloring vector.  Then $(G',v')$ may be obtained from $(G,v)$ by a finite number of the following operations:
$$\Lambda_{1}: \; (G,v)\mapsto (PGP^{T}, Pv \; \modp)$$
and
$$\Lambda_{2}: \; (G,v)\mapsto (G'',v'')$$
where $P$ is an invertible, unimodular, integer matrix and 
$$G''=\left(
\begin{array}{c c c c c}
	\; & \; & \; & \ast & 0 \\
	\; & G & \; & \vdots & \vdots \\
	\; & \; & \; & \ast & 0 \\
	\ast & \cdots & \ast & 0 & 1 \\
	0 & \cdots & 0 & 1 & 0
\end{array}
\right)$$ 
and $v''=\left( \begin{array}{c} v \\ 0 \\ 0 \end{array} \right)$.  A straightforward calculation shows that $cu$ is unchanged by either of the $\Lambda$-moves.

The effect on $(G,v)$ when we add a single twisted band is 
$$G''
= \left(
\begin{array}{c c c c }
	\; & \; & \; & 0  \\
	\; & G & \; & \vdots  \\
	\; & \; & \; & 0  \\
	0 & \cdots & 0 & \pm 1
\end{array}\right) $$
and $v''=\left(\begin{array}{c}v \\ 0 \end{array}\right)$ (see Figure \ref{fig:invariants8}).

\begin{figure}
	\centering
		\includegraphics[scale = .4]{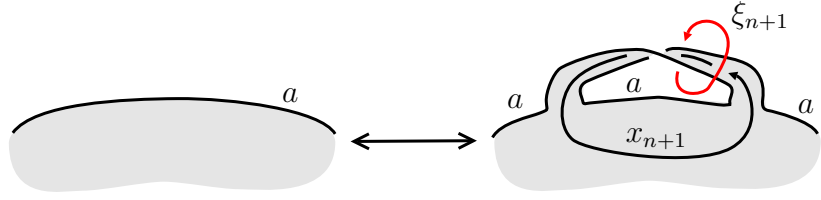}
	\put(-35,65){$\xi_{n+1}$}
	\put(-10,31){$a$}
	\put(-205,35){$a$}
	\put(-75,39){$a$}
	\put(-120,33){$a$}
	\put(-75,20){$x_{n+1}$}
	\caption{Coloring resulting from addition of a twisted band.}
	\label{fig:invariants8}
\end{figure}

Thus, the colored untying invariant defined by the Goeritz matrix is unchanged by any of the moves.
\end{proof}

\begin{figure}[hbtp]


  \vspace{9pt}

  \centerline{\hbox{ \hspace{0.50in}
    \epsfxsize=4.0in
    \epsffile{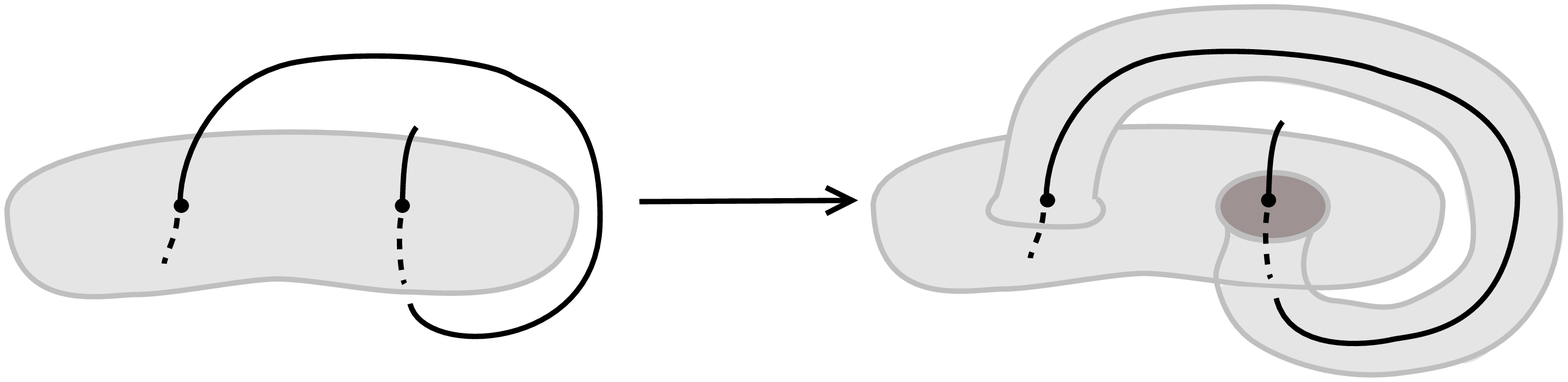}
    	\put(-242,4){\large$F$}
			\put(-255,57){\large$L$}
    }
  }

  \vspace{9pt}

  \caption{``Tubing off'' intersections with the same sign.}
  \label{fig:invariants9-10}

\end{figure}

The next lemma will be used exclusively in the proof of Theorem \ref{thm:cuInv2} below.
\begin{lem}
If $L\subset S^{3}-K$ is a link so that its homotopy class $[L]$  is in $ker(\rho)$ then $L \subset S^{3}-F$ for some spanning surface $F$ for $K$.  Notice that we do not need to assume that $L$ is an unlink. \label{lem:cuDisjoint}
\end{lem}
\begin{proof}
From Diagram \ref{ColoringComDiagram} we have seen that if $[L]$ is in the kernel of $\rho$ then $lk(L,K)\equiv 0 \; (mod \; 2)$.  Then $L$ intersects $F$ an even number of times however two adjacent (innermost) intersections can have opposite or the same sign.  If they have opposite sign then we may resolve them by ``tubing off'' these intersections with a tube which does not change the orientability of the surface.  Otherwise we may resolve the intersections with a ``non-orientable tube'' as in Figure \ref{fig:invariants9-10}.  The resulting spanning surface is $S$-equivalent (in the non-orientable sense of $S$-equivalence) to $F$ and has reduced the number of intersections with $L$.
\end{proof}

\subsubsection{Surgery equivalence}

First we will prove Proposition \ref{prop:cuGoer}, then we will show, via the Goeritz definition, that the colored untying invariant is an invariant of $\pm 1$-framed surgery in the kernel of $\rho$.

\begin{proof}\textbf{\textit{Proposition \ref{prop:cuGoer}}.}\\
The authors would like to thank Pat Gilmer for suggesting this method of proof.

We wish to relate $cu(K,\rho)=a \cup \beta^{1}(a)\in \Z_{p}$ to $cu(K,\rho)'=\frac{v^{T}Gv}{p}\modp$.  We will show that the ``bockstein definition'' $cu(K,\rho)$ is given by the linking pairing on $H^{1}(\tilde{X};\Q/\Z)$ where $\tilde{X}$ is the double-branched cover along $K$ of the $3$-sphere.  On the other hand the Goeritz matrix gives an equivalent linking pairing on $Hom(H_{1}(\tilde{X});\Q/\Z)$.  Moreover, given a presentation of the first homology of the double-branched cover, the two pairings give the same element of $\Z_{p}\subset \Q/\Z$.

Consider the following commutative diagram consisting of coefficient groups.
\begin{equation*}
\xymatrix{
	0 \ar[r] & \Z \ar[r]^{\times p} \ar[d]^{=} & \Z \ar[r]^{\modp} \ar[d]^{\times 1/p} & \Z_{p} \ar[r] \ar[d]^{j} & 0 \\
	0 \ar[r] & \Z \ar[r] & \Q \ar[r] & \Q/\Z \ar[r] & 0 	   
}
\end{equation*}
where $j$ is the natural inclusion of $\Z_{p}$ into $\Q/\Z$, more precisely $\Z_{p}\cong (1/p)\Z/\Z \subset \Q/\Z$.  In particular, if $\hat{a}$ is the element of $H^{1}(\tilde{X},\Q/\Z)$ corresponding to $a\in H^{1}(\tilde{X};\Z_{p})$ from the bockstein definition of the colored untying invariant then $\hat{a}$ is determined by the vector $\hat{v}=\frac{v}{p}$ with respect to a choice of basis for $H_{1}(\tilde{X};\Z)$.  That is $v$ is the coloring vector which describes where the ``coloring'' $\rho'$ sends a generating set of $H_{1}(\tilde{X};\Z)$.  

Under the isomorphisms
\begin{equation}\label{eq:LinkingForm}
	\xymatrix{
	Hom(H_{1}(\tilde{X}),\Q/\Z) \ar[r]^>>>{\Gamma}_>>>{\cong} & H^{1}(\tilde{X};\Q/\Z) \ar[r]^>>>{\beta^{1}}_>>>{\cong} & H^{2}(\tilde{X};\Z) \ar[r]_>>>{\cong} & H_{1}(\tilde{X};\Z)
	}
\end{equation}
arising from the universal coefficient theorem, the definition of the bockstein homomorphism $\beta^{1}$, and Poincare' duality there is a correspondence between the bilinear pairing on $H^{1}(\tilde{X};\Q/\Z)$ defined by $(a,b)\mapsto \left[\tilde{X}\right]\cap (a\cup \beta^{1}(b))$ and the linking form on $H_{1}(\tilde{X};\Z)$.  Here $[M]\in H_{3}(\tilde{X})$ denotes the fundamental class of the $3$-manifold.  Furthermore, under the isomorphism $\Gamma$, the pairing corresponds to the form $\lambda$ given in \cite[page 8]{Gil} on $Hom(H_{1}(\tilde{X}), \Q/\Z)$ relative to the generators $\left\{x_{i}\right\}$ for $H_{1}(F)$, for some spanning surface $F$, and their duals $\left\{\xi_{i}\right\}$ which generate $H_{1}(S^{3}-F)$.  Now by \cite{Gor-Li}, this matrix is the Goeritz matrix $G$.
Thus
\begin{eqnarray*}
p\cdot \lambda(\Gamma^{-1}(\hat{a}),\Gamma^{-1}(\hat{a})) &=& \hat{a}\cup \beta(\hat{a})\\
																										 &=& cu(K,\rho)\in \Z_{p}\subset \Q/\Z.  
\end{eqnarray*}
And so $$\frac{cu(K,\rho)}{p}=\frac{v^{t}Gv}{p^{2}}=\frac{cu(K\rho)'}{p}$$ as desired.
\end{proof}



We will now show that the colored untying invariant is a surgery equivalence invariant for $p$-colored knots.  Note that Moskovich gives an alternate algebraic proof in \cite{Mos}.

\begin{thm}\label{thm:cuInv2}
The colored untying invariant is invariant under $\pm 1$-framed surgery in the kernel of $\rho$. 
\end{thm}

\begin{proof}
From Proposition \ref{prop:cuGoer} we may assume that $cu(K,\rho)=\frac{v^{T}\cdot G \cdot v}{p} \; \modp$ for some coloring vector $v=(v_{1},\ldots, v_{n})^{T}$ and Goeritz matrix $G$ corresponding to a spanning surface for $K$.  Let $[L]$ be in the kernel of the coloring for $K$ represented by an unlink $L$ in the complement of the knot.  Lemmas \ref{lem:cuSurf}, and \ref{lem:cuDisjoint} imply that the spanning surface may be chosen so that $L$ is disjoint from the spanning surface.  Furthermore, let $K$ be in \textit{disk-band form} (see \cite[Chapter 8]{Bur-Zi}).

Under these conditions, $\pm 1$-surgery along one component of $L$ adds a single full twist in $k$ parallel bands of $K$ corresponding to generators (after renumbering perhaps) $x_{1}, \ldots, x_{k}$ for $H_{1}(F)$ with $v_{1}+ \cdots + v_{k}\equiv 0 \; (mod \; p)$.  Then the pair $(G,v)$ changes as follows:
$$G \mapsto G + \left[ \begin{array}{c c} N & 0 \\ 0 & 0 \end{array}\right]=G' \; and \;  v \mapsto v$$
where $N$ is a $k \times k$ matrix whose entries are all $2$.  Thus,
\begin{eqnarray*}
v^{T} G' v &=& p \cdot cu(K,\rho) + v^{T} \left[ \begin{array}{c c} N & 0 \\ 0 & 0 \end{array}\right] v \\
											&=& p\cdot cu(K,\rho) + (v_{1} \; \cdots v_{k})\left[ \begin{array}{c c} N & 0 \\ 0 & 0 \end{array}\right] \left( \begin{array}{c} v_{1} \\ \vdots \\ v_{k} \end{array} \right) \\
											&=& p\cdot cu(K,\rho) + 2(v_{1}+ \cdots +v_{k})^{2} \\
											& \equiv & p\cdot cu(K, \rho) \; (mod \; p^{2})
\end{eqnarray*} 
and so the colored untying invariant is unchanged by $\pm 1$-surgery along $L$.
\end{proof}

We will now show by explicit example that $cu$ is non-trivial for all $p$.  We will also show that there are at least $p$ surgery classes of $p$-colored knots and that connected sums of $(p,2)$-torus knots give a representative of each of these $p$ classes.  Note that once again an alternate proof of these results is in \cite{Mos}.

\subsubsection{Examples} \label{sec:cuEx}

Since we may pick any spanning surface for the knot regardless of orientation, we shall always use the spanning surface corresponding to a checkerboard coloring for a diagram for $K$.

\begin{ex} \textbf{$7$-colorable knots of genus $1$ with at most $12$ crossings.}

\noindent From the table of knots given by \textit{KnotInfo} \cite{Knot}, the only $7$-colorable knots of genus $1$ with at most $12$ crossings are $5_{2}$, $7_{1}$, $11_{n141}$, and $12_{a0803}$.  We will show that the colors of two arcs at any crossing in the diagrams given in Figure \ref{fig:invariants12-15} determine the coloring as well as the colored untying invariants.  Note that Figure \ref{fig:invariants12-15} $(a)$ shows the coloring which is forced by the choice of $a$ and $b$ in $\Zp$ as well as the choice of generators $\{ x_{i} \}$ and $\{ \xi_{i} \}$ for $H_{1}(F;\Z)$ and $H_{1}(S^{3}-F;\Z)$ respectively.  However, in $(b)$-$(d)$, the redundant labels are omitted.  The infinite region is labeled by $\ast$ and the other white regions are understood to be numbered to coincide with the numbering of the $\xi$'s.
   
\begin{figure}[hbtp]
  \vspace{9pt}

  \centerline{\hbox{ \hspace{0.0in} 
    \epsfxsize=2.0in
    \epsffile{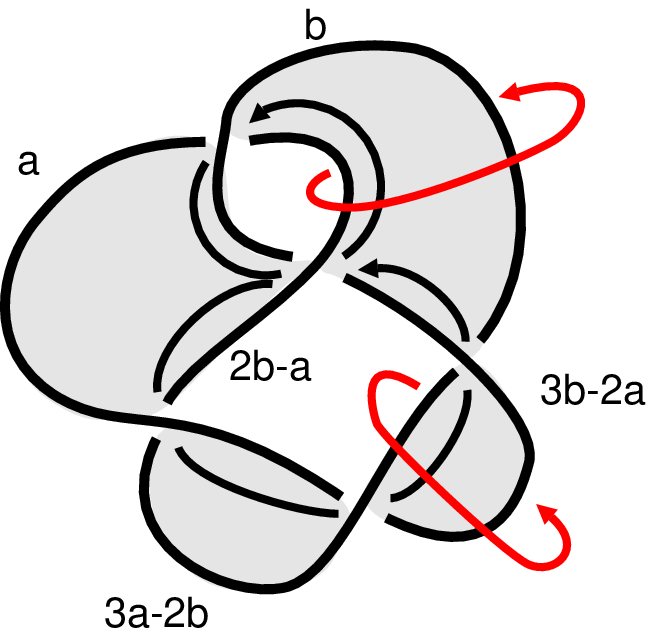}
    	\put(-22,20){$\xi_{2}$}
			\put(-20,110){$\xi_{1}$}
			\put(-92,20){$x_{2}$}
			\put(-120,85){$x_{1}$}
			\put(-150,120){\Huge$\ast$}
    \hspace{0.30in}
    \epsfxsize=2.0in
    \epsffile{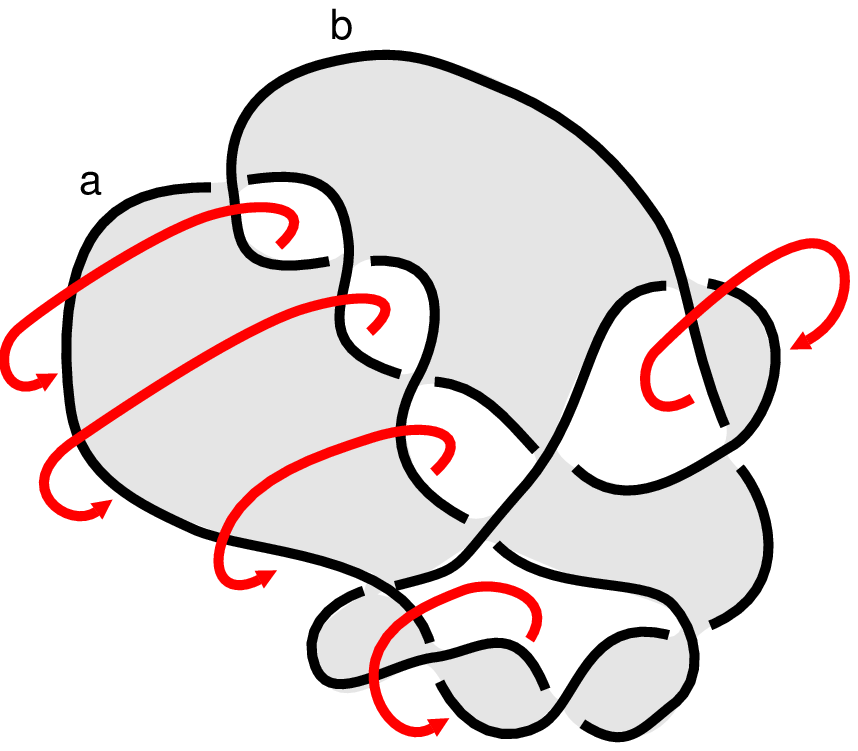}
    	\put(2,95){$\xi_{4}$}
    	\put(-90,10){$\xi_{5}$}
    	\put(-157,80){$\xi_{1}$}
    	\put(-152,55){$\xi_{2}$}
    	\put(-120,40){$\xi_{3}$}
    	\put(-150,120){\Huge$\ast$}
    }
  }

  \vspace{9pt}
  \hbox{\hspace{1.35in} (a) \hspace{2.10in} (b)} 
  \vspace{9pt}

  \centerline{\hbox{ \hspace{0.50in}
    \epsfxsize=2.0in
    \epsffile{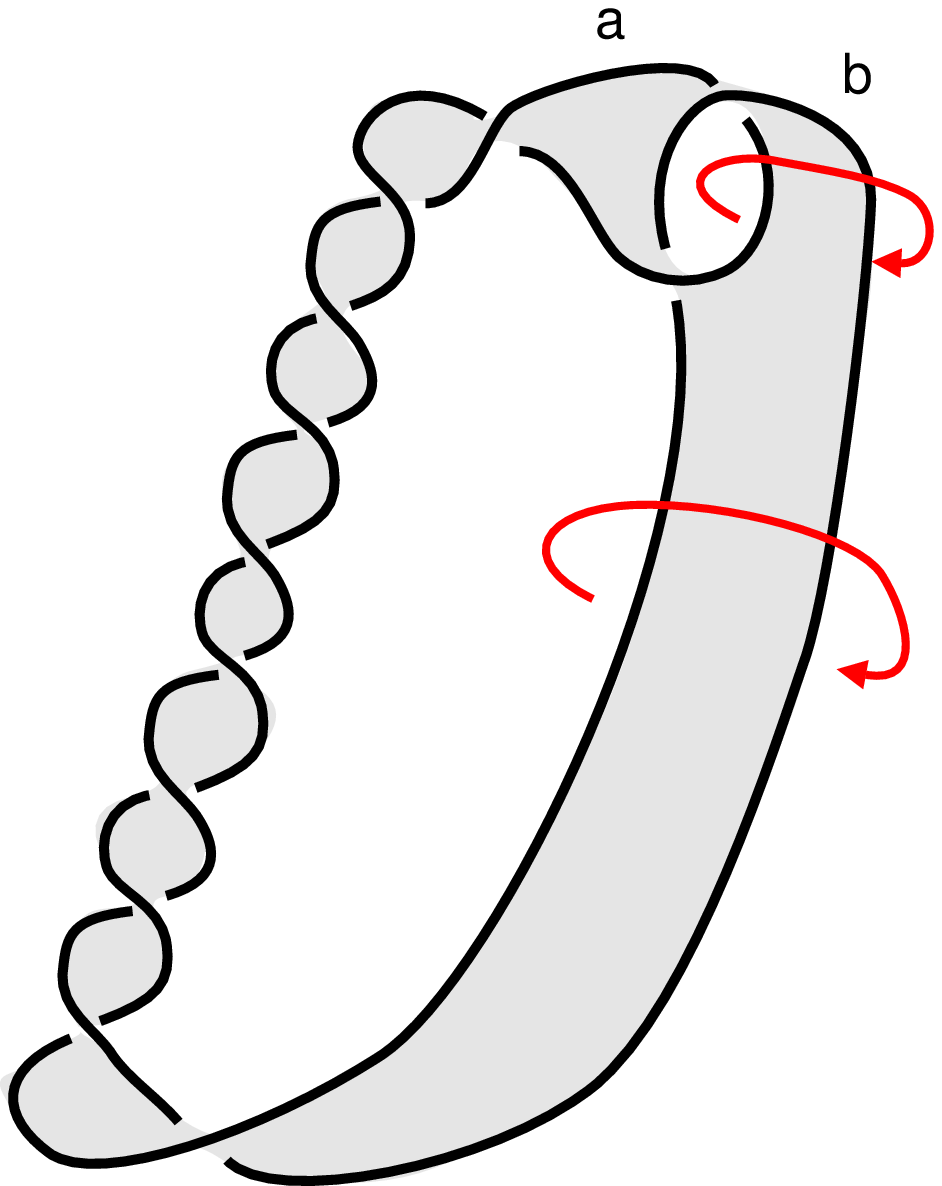}
    	\put(2,140){$\xi_{2}$}
    	\put(0,85){$\xi_{1}$}
    	\put(-150,120){\Huge$\ast$}
    \hspace{0.30in}
    \epsfxsize=2.0in
    \epsffile{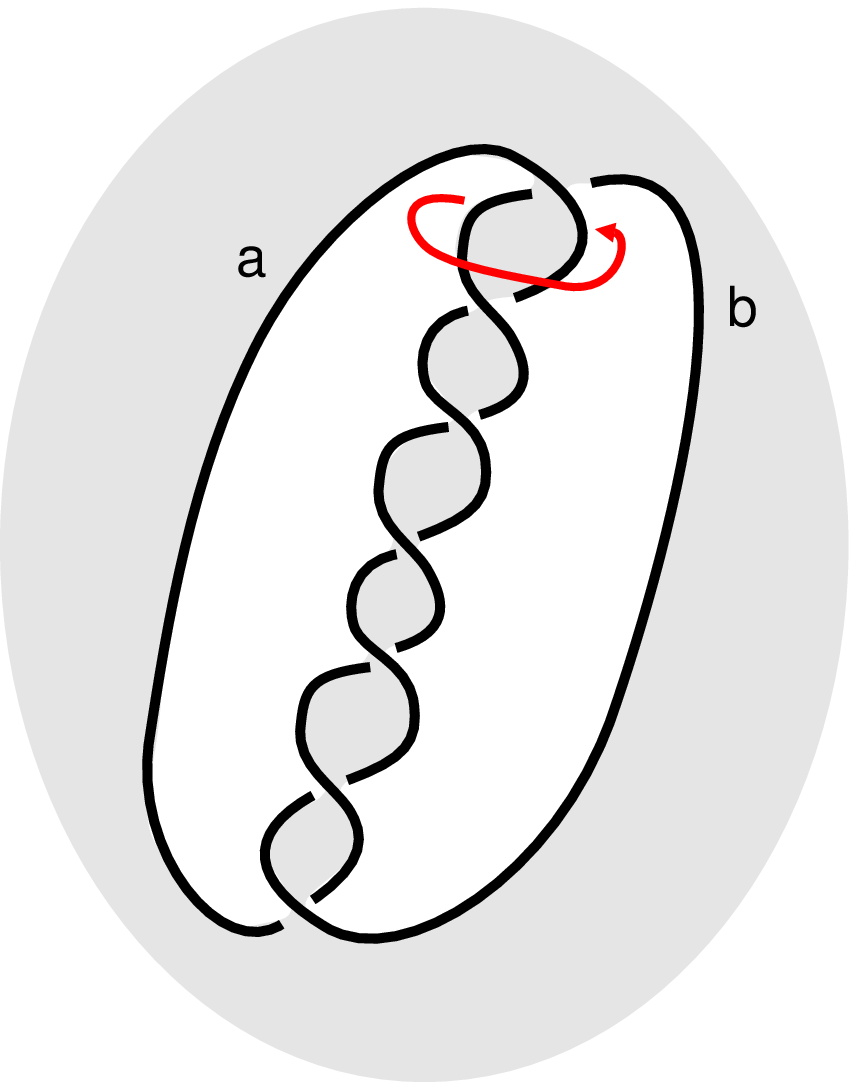}
    	\put(-38,130){$\xi_{1}$}
    	\put(-100,90){\Huge$\ast$}
    }
  }

  \vspace{9pt}
  \hbox{\hspace{1.35in} (c) \hspace{2.10in} (d)} 
  \vspace{9pt}

  \caption{The (a) $5_{2}$, (b) $11_{n141}$, (c) $12_{a0803}$, and (d) $7_{1}$ knots.}
  \label{fig:invariants12-15}

\end{figure}

\begin{prop} \label{ex:cu7Colorable}
The colored untying invariants for the $7$-colorable knots $5_{2}$, $11_{n141}$, $12_{a0803}$, and $7_{1}$ are non-zero multiples of squares for any non-trivial coloring.  In particular, there are three distinct values of $cu$, one for each square modulo $7$, for each of the four knots depending on the coloring class.
\end{prop}

\begin{proof}
First we must pick a white region in a checkerboard coloring for the diagram to be the so-called \textit{infinite region}.  If $F$ is the spanning surface described by the black regions of the checkerboard coloring, then a basis for $H_{1}(F)$ is represented by loops $\left\{ x_{1},\dots , x_{n}\right\}$ which are parallel to the boundary of each white region excluding the infinite region.  Then, the coloring vector is $$v=(\overline{\rho}(\xi_{1}),\ldots,\overline{\rho}(\xi_{n}))^{T}$$ where $\overline{\rho}:H_{1}(S^{3}-F)\rightarrow \Zp$ is the map at the top of Diagram \ref{ColoringComDiagram}, and $q:\Zp\rightarrow \Z$ is the forgetful map as in Proposition \ref{prop:cuGoer}.  Here $\left\{ \xi_{i}\right\}$ is a basis for $H_{1}(S^{3}-F)$ represented by loops in the complement of the surface that pass through the infinite region and the $i$th white region exactly once each so that $lk(x_{i}, \xi_{i})=\delta_{ij}$.

Then the Goeritz matrices in question are:
	$$G(5_{2})=\left(\begin{array}{c c} -2 & 1 \\ 1 & -4 \end{array} \right),$$
	$$G(11_{n141})=
	\left(\begin{array}{c c c c c}  -2 & 1 & 0 & 0 & 0 \\
	 																1 & -2 & 1 & 0 & 0 \\
	 																0 & 1 & -1 & -1 & 1 \\
	 																0 & 0 & -1 & 3 & 0 \\
	 																0 & 0 & 1 & 0 & -5	\end{array} \right),$$
	$$G(12_{a0803})=
	\left(\begin{array}{c c} -11 & 1 \\ 1 & -2 \end{array} \right),$$ and
	$$G(7_{1})=(-7).$$  And so the colored untying invariants are: $cu(5_{2})=5(b-a)^{2}$, $cu(11_{n141})=5(b-a)^{2}$, $cu(12_{a0803})=(b-a)^{2}$, and $cu(7_{1})=6(b-a)^{2}$ where each is understood to be modulo $7$.
\end{proof}
\end{ex}

Notice that the above construction for $cu(7_{1})$ easily generalizes for all odd primes $p$.

\begin{ex}\label{ex:cu(p,2)} \textbf{The $(p,2)$-torus knots for any $p$.}

\noindent The $7_{1}$ knot is also known as the $(7,2)$-torus knot.  As an extension of the construction used to calculate $cu(7_{1})$, Figure \ref{fig:invariants16} gives the general result.  Note that the $p$ in the figure denotes $p$ positive half twists.

\begin{figure}
	\centering
		\includegraphics[scale = .8]{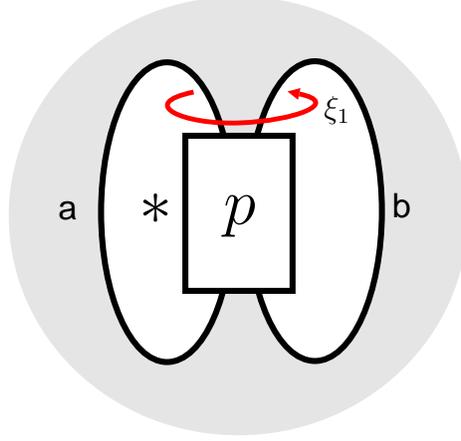}
	\put(-125,80){\huge$\ast$}
	\put(-55,120){$\xi_{1}$}
	\put(-93,80){\huge$p$}
	\caption{The $(p,2)$-torus knot}
	\label{fig:invariants16}
\end{figure}

\noindent So $cu((p,2))=-(b-a)^{2}$ which implies that there is one colored untying class for each square modulo $p$ for the $(p,2)$-torus knot.
\end{ex}

We will now show that the colored untying invariant is additive under the operation of the connected sum of $p$-colored knots.  As an immediate corollary of this we see that the connected sum of $k \; (p,2)$-torus knots for $k=1,\ldots,p$, with the appropriate choices of colorings, give a complete list of representatives of the colored untying invariant classes.

\begin{figure}
	\centering
		\includegraphics[scale = .5]{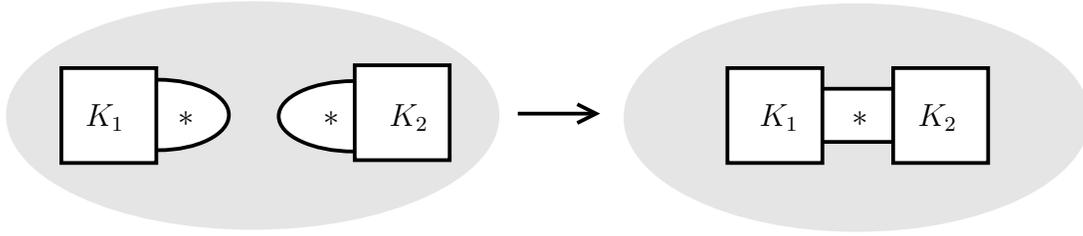}
	\put(-380,40){$K_{1}$}
	\put(-265,40){$K_{2}$}
	\put(-125,40){$K_{1}$}
	\put(-65,40){$K_{2}$}
	\put(-90,40){$\ast$}
	\put(-345,40){$\ast$}
	\put(-290,40){$\ast$}
	\caption{Checkerboard coloring for a connected sum}
	\label{fig:invariants11}
\end{figure}

Note that, as Figure \ref{fig:invariants11} suggests, if we pick an appropriate checkerboard coloring the proof of the following proposition is clear.  Namely, we wish to pick the infinite regions for the checkerboard colorings for the summands so that the checkerboard coloring for the connected sum is determined.
\begin{prop}
The colored untying invariant is additive under the operation of the connected sum of $p$-colored knot.
\end{prop}


We have shown that $cu(K,\rho)$ is a non-trivial, additive, surgery equivalence invariant of $p$-colored knots.  We will now define the \textit{bordism invariants} which exhibit the same properties.  They are all, in fact, the same invariant.  We used the Goeritz definition of the colored untying invariant to establish a lower bound on the number of surgery equivalence classes.  To obtain an upper bound we will need a definition of $cu$ in the context of bordism theory.

\subsection{The bordism invariants} \label{sec:BordismInvariants}

Once again, let $(K,\rho)$ be a $p$-colored knot, and let $\tilde{X}_{0}$ and $X_{0}$ be the manifolds obtained by performing $0$-framed surgery along $K$ to the manifolds $\tilde{X}$ (the $2$-fold brached cover of $S^{3}$) and $S^{3}$ itself respectively.  If we have a map $f:M^{3}\rightarrow K(G)$ where $K(G)$ denotes the Eilenberg-Maclane space $K(G,1)$ then the image of the fundamental class under the induced map $f_{\ast}: H_{3}(M; \Z)\rightarrow H_{3}(K(G); \Z)$ is an invariant of the $3$-manifold $M$.  The construction is exactly the same as the invariants described by T. Cochran, A. Gerges, and K. Orr in \cite{CGO}.  We will divide the bordism invariants into two categories \textit{closed} and \textit{relative}.

\subsubsection{The closed bordism invariants}\label{sec:ClosedBordInv}

As mentioned earlier in the ``Preliminaries'' $$H_{3}(K(G); \Z)\cong \Omega_{3}(G)$$ and it is in this context that the \textit{bordism}  \textit{invariants} arise.  So to define $\omega_{2}$ and $\omega_{0}$ we must find maps from $\tilde{X_{0}}$ and $X_{0}$ to Eilenberg-Maclane spaces over the appropriate groups.

We wish to have maps which arise naturally from the coloring $\rho$.  Recall that the second derived group of $G$, denoted $G^{(2)}$, is defined to be the commutator subgroup of the commutator subgroup of $G$.  That is $G^{(2)}=\left[G_{2},G_{2}\right]$ where $G_{2}=\left[G,G\right]$.  Since a preferred longitude of the knot $K$ is in the second derived group of $\pi_{1}(S^{3}-K)$ it must be mapped trivially by $\rho$.  Hence the map $\rho'$ from Diagram \ref{ColoringComDiagram} factors through
\begin{equation*}
\xymatrix{
	\pi_{1}(\tilde{X}_{0}) \ar[ddrr] & \; & \; \\
	\pi_{1}(\tilde{X}-\tilde{K}) \ar[u] & \; & \; \\ 
	\pi_{1}(\tilde{X})\ar[u] \ar[rr]^{\rho'} & \; & \Zp	  
}
\end{equation*}
\noindent which establishes the existence of a map $\tilde{f}:\tilde{X}_{0}\rightarrow K(\Zp)$ as desired.  Likewise, and perhaps even easier to see, we have that $\rho$ factors through
\begin{equation*}
\xymatrix{
	\pi_{1}(X_{0}) \ar[drr]& \; & \; \\
	\pi_{1}(S^{3}-K)\ar[u] \ar[rr]^{\rho} & \; & D_{2p}	  
}
\end{equation*}
\noindent which gives a natural map $$f:X_{0}\rightarrow K\left(\frac{\pi_{1}(S^{3}-K)}{ker \; \rho}\cong D_{2p},1\right).$$  We will show that the induced maps on homology of $\tilde{f}$ and $f$ define invariants of not only the $3$-manifolds $\tilde{X_{0}}$ and $X_{0}$ but they are also surgery equivalence invariants for the $p$-colored knot $(K,\rho)$.

\begin{definition}
Suppose $\tilde{f}:\tilde{X_{0}}\rightarrow K(\Zp)$ and $f:X_{0}\rightarrow K(D_{2p})$ are the maps obtained via the coloring $\rho$ as above.  Then define the \textit{closed bordism invariants} to be
$$\omega_{2}(K,\rho):=\tilde{f}_{\ast}([\tilde{X_{0}}])\in H_{3}(\Zp; \Z)$$
and
$$\omega_{0}(K,\rho):= f_{\ast}([X_{0}])\in H_{3}(D_{2p}; \Z)$$
\noindent where $[M]\in H_{3}(M; \Z)$ denotes the fundamental class of $M$.  
\end{definition} 

Notice that the invariants depend only on the bordism classes of the (closed) $3$-manifolds over $\Zp$ and $D_{2p}$ respectively which is the motivation for the names.  It is also clear that $\tilde{X}$ and $\tilde{X_{0}}$ are in the same bordism class over $\Zp$.  The bordism is constructed from $\tilde{X}\times [0,1]$ by attaching a $2$-handle along the lift of the prefered longitude.

The final bordism invariant, denoted simply by $\omega$, arises from the manifold $M=(S^{3}-K)$ which is not closed so it will be defined separately.  Also note that since $\pm 1$-framed surgery along links in the kernel of the coloring $\rho$ defines a bordism between the resulting manifolds then $\omega_{2}$ and $\omega_{0}$ are actually surgery equivalence invariants.  The bordism is obtained by attaching a $2$-handle along each component of the surgery link to $M\times [0,1]$ (for $M=\tilde{X_{0}}, X_{0}$).

\subsubsection{The relative bordism invariant}

Recall the definition of a \textit{based $p$-colored knot} which is a $p$-colored knot with a chosen meridian $m$ so that $\rho(m)=ts^{0}$.  That is, if the coloring $\rho$ is defined by a labeling of a diagram for $K$ then the arc corresponding to $m$ would have the label $0$.  We may assume this because $p$-colored knots are only defined up to an inner automorphism of the dihedral group.  This allows, in particular, for any chosen arc to have the label $0$.  We will now define the last of the three bordism invariants.

\begin{definition}
Let $(K,\rho,m)$ be a based $p$-colored knot.  If $K(\mathbb{Z}_{2})$ is the subspace of $K(D_{2p})$ corresponding to the image of $m$ under the coloring, then define
$$\omega(K,\rho):=\rho([M,\partial M])\in H_{3}(K(D_{2p}),K(\mathbb{Z}_{2});\Z)$$
\noindent where $[M,\partial M]$ denotes the fundamental class of $M=(S^{3}-K)$ relative to the boundary and $f:(M,\partial M ) \rightarrow (D_{2p},\mathbb{Z}_{2})$ arises directly from the coloring.
\end{definition}

\noindent Indeed, we may think of $K(\mathbb{Z}_{2})$ as a subspace of $K(D_{2p})$ because we may construct a $K(D_{2p})$ from a $K(\mathbb{Z}_{2})$ by adding $k$-cells, $k=1,2,\dots$, to obtain the correct homotopy groups.  Furthermore, since we can assume that the fundamental group of the boundary torus is generated by the classes represented by the preferred longitude and our chosen meridian $m$, it is clear that $\partial M$ is mapped into the correct subspace.

We will now prove a few special properties of the bordism invariants.

\subsubsection{Properties}

Consider the \textit{bordism long exact sequence of the pair} $(X,A)$
\begin{equation}\label{BordismLES} \cdots \longrightarrow \Omega_{n}(A) \stackrel{i_{\ast}}{\longrightarrow} \Omega_{n}(X) \stackrel{j_{\ast}}{\longrightarrow} \Omega_{n}(X,A)\longrightarrow \Omega_{n-1}(A)\longrightarrow \cdots 
\end{equation}
for $i_{\ast}$ and $j_{\ast}$ induced by inclusion (see Section $5$ of \cite{Co-Fl}).  We will be concerned with the pairs $(X,A)=(K(D_{2p}),K(\Zp))$ and $(X,A)=(K(D_{2p}), K(\Z_{2}))$ which will relate $\omega_{2}$ to $\omega_{0}$, and $\omega_{0}$ to $\omega$ respectively.

In these cases, we may compute the bordism groups using the fact that $\Omega_{n}(K(G,1))\cong H_{n}(G;\Z)$.  The cohomology groups of cyclic groups are well-known and may be computed using a spectral sequence for the \textit{fibration}
$$K(\Z,1)\rightarrow K(\Zp,1) \rightarrow K(\Z, 2)$$ with fiber $K(\Z,1)$ being a circle (see Chapter 9 \cite{Da-Ki}).  The homology groups are then obtained from the cohomology groups by using the Universal Coefficient Theorem.  We have 
\begin{displaymath}
	H_{n}(\Zp) \cong
	\begin{cases}
		\Z & \text{if $n=0$,} \\
		\Zp & \text{if $n$ is odd, and} \\
		0 & \text{if $n>0$ is even}
	\end{cases}
\end{displaymath}
for $p$ any prime number.  

The following proposition follows from a spectral sequence found in \cite{Ad-Mil} and is well-known.  

\begin{prop} \label{prop:HomDihedral}
The homology groups of the dihedral group $D_{2p}$ are as follows
\begin{displaymath}
  H_n(D_{2p}) \cong
  \begin{cases}
    \Z_2 & \text{if $n \equiv 1 \pmod{4}$,}\\
    \Z_2 \oplus \Zp & \text{if $n \equiv 3 \pmod{4}$} \\
    0 & \text{otherwise}
  \end{cases}
\end{displaymath}%
if $p$ is an odd prime.
\end{prop}   
\noindent So the closed bordism invariants $\omega_{2}$ and $\omega_{0}$ may be thought of as elements of $\Z_{p}$ and $\Z_{2p}$ respectively.

We will use the bordism long exact sequence \ref{BordismLES} to determine the group in which the relative bordism invariant $\omega$ resides.  Consider \begin{equation*} \cdots \longrightarrow \Omega_{3}(K(\Z_{2})) \stackrel{i_{\ast}}{\longrightarrow} \Omega_{3}(K(D_{2p})) \stackrel{j_{\ast}}{\longrightarrow} \Omega_{3}(K(D_{2p}),K(\Z_{2}))\longrightarrow \Omega_{2}(K(\Z_{2}))\longrightarrow \cdots 
\end{equation*}
\noindent where the Eilenberg-Maclane space $K(\Z_{2})$ is the subspace of $K(D_{2p})$ arising from the subgroup $\Z_{2}\cong \left\langle t \right\rangle \in D_{2p}$.  In this case $i_{\ast}$ is injective since any \textit{singular manifold} $(M, \varphi)$ that is null-bordant over $D_{2p}$ is null-bordant over $\Z_{2}$ via the same $4$-manifold.  As $\Omega_{2}(K(\Z_{2}))$ is trivial we have 
\begin{equation*}
	0 \rightarrow \Omega_{3}(K(\Z_{2}))\cong \Z_{2} \hookrightarrow \Omega_{3}(K(D_{2p}))\cong \Z_{2p}\rightarrow \Omega_{3}(K(D_{2p}), K(\Z_{2})))\rightarrow 0
\end{equation*}
\noindent is exact.  In particular $\Omega_{3}(K(D_{2p}), K(\Z_{2})) \cong \Zp$.  So the relative bordism invariant $\omega$ may be regarded as an element of $\Zp$.  We will later show, in the proof of Theorem \ref{thm:Main}, that the closed bordism invariant $\omega_{0}\in \Omega_{3}(K(D_{2p}))\cong \Z_{2p}$ only takes values in the $\Zp$ part of $\Z_{2p}$ which will establish an equivalence between all three bordism invariants.

We have already seen that the colored untying invariant is additive under the operation of the connected sum of $p$-colored knots.  The same is true for the bordism invariants.  Of course, once we have established the equivalence of all the invariants, then the additivity of $cu$ is enough to show this and a direct proof will be omitted.  

As a corollary to the additivity of the closed bordism invariants we see that if $\omega_{0}\in \Omega_{3}(K(D_{2p}))\cong \Z_{2}\oplus \Zp$ is a $p$-valued invariant of $p$-colored knots, then $\omega_{0}(K,\rho)=(0,n)$ for any colored knot $(K,\rho)$.  This is because if $\omega_{0}$ is $p$-valued then every $p$-colored knot must have the same value in the first coordinate of $\Z_{2}\oplus \Zp \cong \Z_{2p}$.  And since $\omega_{0}(K\# K, \rho \# \rho)=(0,2n)$, the first coordinate value must be $0$.  We will show that $\omega_{2}(K,\rho)=2\omega_{0}(K,\rho)$ which will establish an equivalence between $\omega_{0}$ and $\omega_{2}$ once we show that $\omega_{0}$ is $p$-valued.

\section{Proof of equivalence}

We will now show that all of the $p$-colored knot invariants defined above are the same.

\subsection{Equivalence of $cu$ and $\omega_{2}$} \label{sec:cuOmega2}

\begin{prop} \label{prop:cuAndOmega2}
The colored untying invariant $cu(K,\rho)$ is equivalent to the ($2$-fold branched cover) closed bordism invariant $\omega_{2}(K,\rho)$ for any $p$-colored knot $(K,\rho)$.
\end{prop}

\begin{proof}

Again, denote by $\tilde{X}$ the $2$-fold branched cover of $S^{3}$.  Then by the commutativity of Diagram \ref{ColoringComDiagram}, there is a map $\tilde{f}:\tilde{X}\rightarrow K(\Zp,1)$ which corresponds to the coloring $\rho$.  Let $\beta^{1}: H^{1}(\tilde{X};\Zp)\rightarrow H^{2}(\tilde{X};\Z)$ be a Bockstein homomorphism associated with the coffecient sequence $$0\longrightarrow \Z\stackrel{\times p}{\longrightarrow} \Z \stackrel{mod \; p}{\longrightarrow}\Zp \longrightarrow 0.$$  Recall that if $a\in H^{1}(\tilde{X};\Zp)$ is the cohomology class corresponding to $\rho':H_{1}(\tilde{X};\Z)\rightarrow \Zp$ then $$cu(K,\rho)=a \cup \beta^{1}(a) \in H^{3}(\tilde{X};\Zp)\cong \Zp$$ by Moskovich's definition of the colored untying invariant.  Notice that the identification of $cu(K,\rho)$ with an element of $\Zp$ is via evaluation on the fundamental class.

Consider the maps $\tilde{X}\stackrel{\tilde{f}}{\rightarrow} K(\Zp) \stackrel{id}{\rightarrow} K(\Zp)$.  Thus we have the following commutative diagram:
\begin{equation*}
\xymatrix{
	H_{1}(\tilde{X}) \ar[drr]^{\rho'}\ar[d]_{\tilde{f}_{\ast,1}} & \; & \; \\
	H_{1}(\tilde{X}) \ar[rr]_{i} & \; & \Zp= H_{1}(K(\Zp)) \;	  
}
\label{rho'ComDiagram}
\end{equation*}
\noindent where $i:H_{1}(K(\Zp))\rightarrow \Zp$ corresponds to the cohomology class in $H^{1}(K(\Zp);\Zp)$ induced by the identity $id:K(\Zp)\rightarrow K(\Zp)$.  Notice that $a$ corresponds with the homomorphism $\rho'$ by construction, while $\rho'$ corresponds with the cohomology class $\tilde{f}^{\ast,1}\in H^{1}(\tilde{X}; \Zp)$.  The correspondence of $\tilde{f}^{\ast,1}$ and $a$ is exactly $\tilde{f}^{\ast,1}(i)=a$.

Then, by the properties of cup products we have $$\tilde{f}^{\ast,3}(i \cup \beta^{1}(i))=a \cup \beta^{1}(a)$$ which gives the element of $\Zp$ $\left[\tilde{X}\right]\cap (a \cup \beta^{1}(a))$.  On the other hand, if we think of $(i \cup \beta^{1}(i))$ as a chosed fixed generator of $H^{3}(K(\Zp);\Zp)$, then this is the same as $$\tilde{f}_{\ast,3}\left(\left[\tilde{X}\right]\right) \cap (i\cup \beta^{1}(i))$$ which is the identification of $\omega_{2}(K,\rho)$ with an element of $\Zp$.  Note that the non-triviality of the colored untying invariant implies that $(i \cup \beta^{1}(i))$ is a generator of $H^{3}(K(\Zp);\Zp)$.

Hence, with these identifications of $H^{3}(\tilde{X};\Zp)$ and $H_{3}(K(\Zp);\Z)$ with $\Zp$, the elements $cu(K,\rho) \in H^{3}(\tilde{X};\Zp)$ and $\omega(K,\rho)\in H_{3}(K(\Zp);\Z)$ are the same as elements of $\Zp$.

\end{proof}

\subsection{Equivalence of the bordism invariants}

To show that the closed bordism invariants $\omega_{0}$ and $\omega_{2}$ are equivalent it suffices to show two facts.  First we must show that $\omega_{2}(K,\rho)$ is roughly speaking ``twice'' $\omega_{0}(K,\rho)$.  Then we must show that $\omega_{0}$ is a $p$-valued invariant.  This, in turn, will show that all of the bordism invariants are equivalent to each other and to the colored untying invariant.

\begin{lem}\label{lem:Omega2vsOmega01}
The closed bordism invariants have the property that $\omega_{2}(K,\rho)=2n$ if $\omega_{0}(K,\rho)=(m,n)\in \Z_{2}\oplus \Zp$.
\end{lem}
\begin{proof}
Recall the bordism long exact sequence \ref{BordismLES} \begin{equation*} \cdots \longrightarrow \Omega_{n}(A) \stackrel{i_{\ast}}{\longrightarrow} \Omega_{n}(X) \stackrel{j_{\ast}}{\longrightarrow} \Omega_{n}(X,A)\longrightarrow \Omega_{n-1}(A)\longrightarrow \cdots 
\end{equation*}
\noindent with $X$ and $A$ the Eilenberg-Maclane spaces over $D_{2p}$, $\Zp$, and $\Z_{2}$ where appropriate.  We have
$$
\xymatrix{
	0 \ar[r] & \Omega_{3}(\Zp)    \ar[r]^{i_{\ast}}     & \Omega_{3}(D_{2p}) \ar[r]        & \Omega_{3}(D_{2p},\Zp)    \ar[r] & 0
} $$
so we must show that $i_{\ast}[\tilde{X_{0}}]=2[X_{0}]$.  But $\tilde{X_{0}}$ is a double cover of $X_{0}$ so the result follows. 
\end{proof}

We will now show that all of the $p$-colored knot invariants give the same information.  In particular, this shows that computation of the bordism invariants may be done by computing the colored untying invariant using the Goeritz matrix.

\begin{thm}\label{thm:Main}

All of the $p$-colored knot invariants are equivalent.
\end{thm}
\begin{proof}
By Propostion \ref{prop:cuAndOmega2} we have that for an appropriate choice of generator for $\Zp$ the elements $\omega_{2}(K,\rho)$ and $cu(K,\rho)$ are equal.  By Lemma \ref{lem:Omega2vsOmega01} above we need only show that $\omega_{0}(K,\rho)$ lies in the $\Zp$ part of $\Z_{2p}$ to show that both of the closed bordism invariants are the same.  The final equivalence between $\omega$ and $\omega_{0}$ will follow from the Bordism Long Exact Sequence.

There is a canonical short exact sequence $$0\rightarrow \Zp=\left\langle s \right\rangle \stackrel{\Phi}{\longrightarrow} D_{2p} = \left\langle s,t \; | \; t^{2}=s^{p}=tsts=1\right\rangle \stackrel{\Psi}{\longrightarrow} \Z_{2} \rightarrow 0$$
where $\Z_{2}$ is the cokernel of the map $\Phi$.  As a result, we may construct the commutative diagram
$$
\xymatrix{
	\; & \; & \Zp \ar[d]^{\Phi} \\
	\pi_{1}(X_{0})\ar[rr]^{\rho} \ar[d]_{\alpha} \ar[drr]^{l}& \; & D_{2p} \ar[d]^{\Psi} \\
	\Z \ar[rr]^{(mod \; 2)} & \; & \Z_{2}  
} $$
\noindent where $\rho$ is the coloring applied to the $0$-surgered manifold, $\alpha$ is abelianization, and $l(x)=lk(x, K) \; (mod \; 2)$.  Hence, we have a commutative diagram of the corresponding spaces
$$
\xymatrix{
	\; & \; & K(\Zp) \ar[d] \\ 
	X_{0} \ar[rr]^{f} \ar[d]_{A} \ar[drr]& \; & K(D_{2p}) \ar[d]^{g} \\
	S^{1}=K(\Z)\ar[rr] & \; & \mathbb{R}P^{\infty}=K(\Z_{2})
} $$
\noindent which induces
$$
\xymatrix{
	\; & \; & \Zp \ar[d] \\ 
	\Z = \left\langle \Lambda \right\rangle \ar[rr]^{f_{\ast}} \ar[d]_{A_{\ast}} \ar[drr]& \; & \Z_{2}\oplus \Zp \ar[d]^{g_{\ast}} \\
	0 \ar[rr] & \; & \Z_{2}
} $$
\noindent on the third homology groups.
From this, we see that $\omega_{0}(K,\rho)=f_{\ast}(\Lambda)=(0,n)\in \Z_{2}\oplus \Zp$ for some $n \in \Zp$ since $A_{\ast}=0$.  Note that $g_{\ast}\neq 0$ since
$$ 
\xymatrix{
	\Z_{2}\ar[rr] \ar[drr]^{id}& \; & D_{2p} \ar[d]^{\Psi} \\
	\; & \; & \Z_{2}
} $$
\noindent commutes.  So the closed bordism invariants are equivalent $p$-valued invariants of $p$-colored knots.  This also implies, in particular, that $\omega_{0}$ and the relative bordism invariant $\omega$ must be the same.

For any based $p$-colored knot, the Bordism Long Exact Sequence gives the exact sequence \begin{equation*} \cdots \longrightarrow \Omega_{3}(K(\Z_{2})) \stackrel{i_{\ast}}{\longrightarrow} \Omega_{3}(K(D_{2p})) \stackrel{j_{\ast}}{\longrightarrow} \Omega_{3}(K(D_{2p}),K(\Z_{2}))\longrightarrow \Omega_{2}(K(\Z_{2}))\longrightarrow \cdots 
\end{equation*} 
that is, we have the short exact sequence
\begin{equation*} 0 \longrightarrow \Z_{2} \stackrel{i_{\ast}}{\longrightarrow} \Z_{2}\oplus \Zp \stackrel{j_{\ast}}{\longrightarrow} \Zp \longrightarrow 0 
\end{equation*} which gives an isomorphism between the $\Zp$-part of $\Omega_{3}(D_{2p})$ and $\Omega_{3}(K(D_{2p}),K(\Z_{2}))$ and the result follows.
\end{proof}

Incidentally, as a corollary to the proof of the Theorem we have the following result.  A detailed proof will not be given here but the result follows from the bordism long exact sequence and the fact that $[\mathbb{R}P^{3},\varphi]\neq 0 \in \Z_{2}\cong \Omega_{3}(K(\Z_{2}))$.
\begin{thm}
The bordism group $\Omega_{3}(D_{2p}) \cong \Z_{2}\oplus \Zp$ is generated by the bordism class represented by the disjoint union of the singular manifolds $(\mathbb{R}P^{3},\varphi)$ and $(X_{0},f)$ where $X_{0}$ is the manifold obtained via $0$-surgery along some prime $p$-colored knot $(K,\rho)$ with non-zero bordism invariant (a $(p,2)$-torus knot for example).  The maps $f$ and $\varphi$ correspond to the coloring $$\rho:\pi_{1}(X_{0})\rightarrow D_{2p},$$ and the inclusion $$\phi:\pi_{1}(\mathbb{R}P^{3})\cong \Z_{2} \rightarrow D_{2p}$$ on the fundamental groups respectively. 
\end{thm}

We have shown that there are at least $p$ surgery equivalence classes of $p$-colored knots, we will now show that twice that is an upperbound on the number of equivalence classes.

\subsection{Main result} \label{sec:SurgEquiv}

We would like to show that the colored untying invariant is a complete invariant for $p$-colored knot surgery type.  This is Moskovich's conjecture, since as we have seen, $cu(K,\rho)$ is $p$-valued.  To show that $cu$ is complete we must show that if $cu(K_{1},\rho_{1})=cu(K_{2},\rho_{2})$ then $(K_{1},\rho_{1})$ and $(K_{2},\rho_{2})$ are surgery equivalent.  The main result of this paper is that this indeed is the case at least half of the time.





\begin{figure}
	\centering
		\includegraphics[scale = .4]{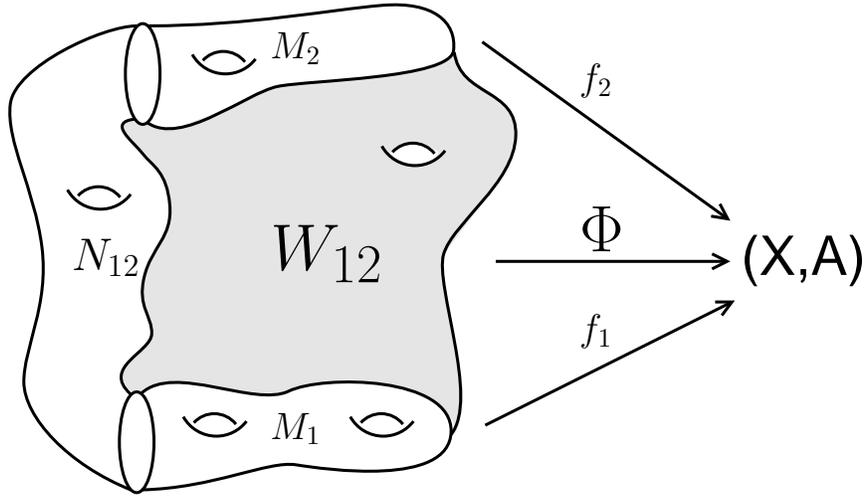}
	\put(-310,85){\Large$N_{12}$}
	\put(-235,85){\huge$W_{12}$}
	\put(-235,165){\large$M_{2}$}
	\put(-235,21){\large$M_{1}$}
	\put(-118,92){\huge$\Phi$}
	\put(-118,58){\large$f_{1}$}
	\put(-118,153){\large$f_{2}$}
	\caption{The ``connecting'' manifold $N_{12}$.}
	\label{fig:surgequiv2}
\end{figure}

Let $P_{a}$ denote the set of all based $p$-colored knots $(K,\rho)$ with $\omega(K,\rho)=a\in \Zp$.  If $(K_{1},\rho_{1})$ and $(K_{2}, \rho_{2})$ are  in the set $P_{a}$ and $M_{i}=S^{3}-K_{i}$ then $(M_{1},\partial M_{1}, f_{1})$ is bordant to $(M_{2},\partial M_{2}, f_{2})$ over $(K(D_{2p},1), K(\Z_{2},1))=(X,A)$ by the definition of the bordism invariant $\omega$.  Here, the $f_{i}:(M_{i},\partial M_{i})\rightarrow (K(D_{2p},1),K(\Z_{2},1))$ are maps which induce the colorings on $\pi_{1}$.  We have the existence of a $4$-manifold $W_{12}$ and a map $$\Phi:(W_{12},\partial W_{12})\rightarrow (K(D_{2p},1),K(\Z_{2},1))$$ so that $\partial W_{12}= (M_{1}\coprod -M_{2})\cup_{\partial N_{12}}N_{12}$ and $\Phi |_{M_{i}}=f_{i}$ as in Figure \ref{fig:surgequiv2}.  The ``connecting'' $3$-manifold in the boundary of the bordism $W_{12}$ between $M_{1}$ and $M_{2}$ is denoted by $N_{12}$.  

Note that the boundary of $N_{12}$ consists of two disjoint copies of the torus $T^{2}$, one for each boundary torus of the $M_{i}$'s.  We would like to show that $N_{12}$ is the product space $T^{2}\times [0,1]$.  We will show that this is necessarily the case at least half of the time.  More precisely, we will construct a map $\eta:P_{a}\times P_{a}\rightarrow \Z_{2}$ that satisfies a certain ``triangle equality'' (Proposition \ref{prop:EtaTriangle} below).  Let $\left[N_{12}\cup_{T^{2}\times \{0,1\}} (T^{2}\times [0,1])\right]=\left[\overline{N_{12}}\right]$ denote the fundamental class of $\overline{N_{12}}$ and let $\Phi_{12}$ be the obvious extension of the map $\Phi |_{N_{12}}:N_{12}\rightarrow K(\Z_{2},1)$ coming from the bordism and $\Phi_{12,\ast}$ denote the induced map on homology.  Define $$\eta(K_{1},K_{2})=\Phi_{12,\ast}(\left[N_{12}\cup_{T^{2}\times \{0,1\}} (T^{2}\times [0,1])\right])$$ which is an element of the bordism group $\Omega_{3}(\Z_{2})\cong \Z_{2}$.  

\begin{prop}
The function $$\eta(K_{1}, K_{2})=0$$ if and only if there is a bordism $(W',\partial W', \Phi')$ between $(S^{3}-K_{1},f_{1})$ and $(S^{3}-K_{2},f_{2})$ with the connecting manifold consisting of the product space $T^{2}\times[0,1]$.
\end{prop}
\begin{proof}
Assume that $\eta(K_{1},K_{2})=0$, note that we must also assume that both knots lie in the set $P_{a}$ in order for the function $\eta$ to make sense.  Then we have bordisms $(W_{0},\Phi_{0})$ over $\Z_{2}$ with boundary $\overline{N_{12}}$ and $(W,\partial W, \Phi)$ over $(D_{2p},\Z_{2})$ with boundary $(S^{3}-K_{1})\cup N_{12}\cup (S^{3}-K_{2})$.  So sufficiency is seen by gluing the bordism $(W_{0},\Phi_{0})$ to the bordism $(W,\partial W,\Phi)$ along the $3$-manifold $N_{ij}$.  The result is a new bordism $(W',\partial W',\Phi')$ over $(D_{2p},\Z_{2})$ defined by 
$$W'=W\cup_{\psi}W_{0}$$
\noindent where $\psi:N_{ij}\rightarrow N_{ij}$ is an orientation reversing diffeomorphism.  The map $\Phi'$ is defined by 
\begin{displaymath}
	\Phi'(x)=
	\begin{cases}
		\Phi_{0}(x) & \text{if $x\in W_{0}$,} \\
		\Phi(x) & \text{if $x\in W$}
	\end{cases}
\end{displaymath} and since the manifolds are glued by a diffeomorphism, it follows that $$\Phi':(W',\partial W')\rightarrow (D_{2p},\Z_{2})$$ is a differentiable map as required.  We have shown that if $\eta(K_{i}, K_{j})$ is trivial then there is a bordism $(W',\partial W',\varphi')$ over $(D_{2p},\Z_{2})$ with $\partial W'=(M_{i}\coprod -M_{j})\cup_{T^{2}\times \{0,1\}} (T^{2}\times [0,1])$.  And so we may assume that $N_{ij}=T^{2}\times [0,1]$ only in the case that $\eta$ is trivial.  Necessity of this condition follows from the fact that $T^{2}\times [0,1]\times [0,1]$ has boundary homeomorphic to $(T^{2}\times[0,1])\cup_{T^{2}\times \{0,1\}}(T^{2}\times [0,1])$.  
\end{proof}

Appealing to the proof of Theorem \ref{thm:MainTheorem} below, we have that the colored surgery untying conjecture in \cite[Conjecture $1$]{Mos} is equivalent to the property that $\eta$ is always trivial.

\begin{cor}
The function $\eta$ vanishes for all pairs of colored knots in $P_{a}$ for all $a\in \Z_{p}$ if and only if there are exactly $p$ surgery equivalence classes of $p$-colored knots.
\end{cor}

We will now show that the map $\eta$ is well-defined and satisfies the ``triangle equality'' property mentioned above. 

\begin{prop} \label{prop:EtaTriangle}
The map $\eta:P_{a}\times P_{a}\rightarrow \Z_{2}$ is well-defined and satisfies 
\begin{equation}\label{eq:EtaTriangle}
\eta(K_{1}, K_{2})=\eta(K_{1},K_{3})+\eta(K_{3},K_{2}) \end{equation}
for any $(K_{3},\rho_{3})\in P_{a}$.
\end{prop}

\begin{proof}


\begin{figure}
	\centering
		\includegraphics[scale = .8]{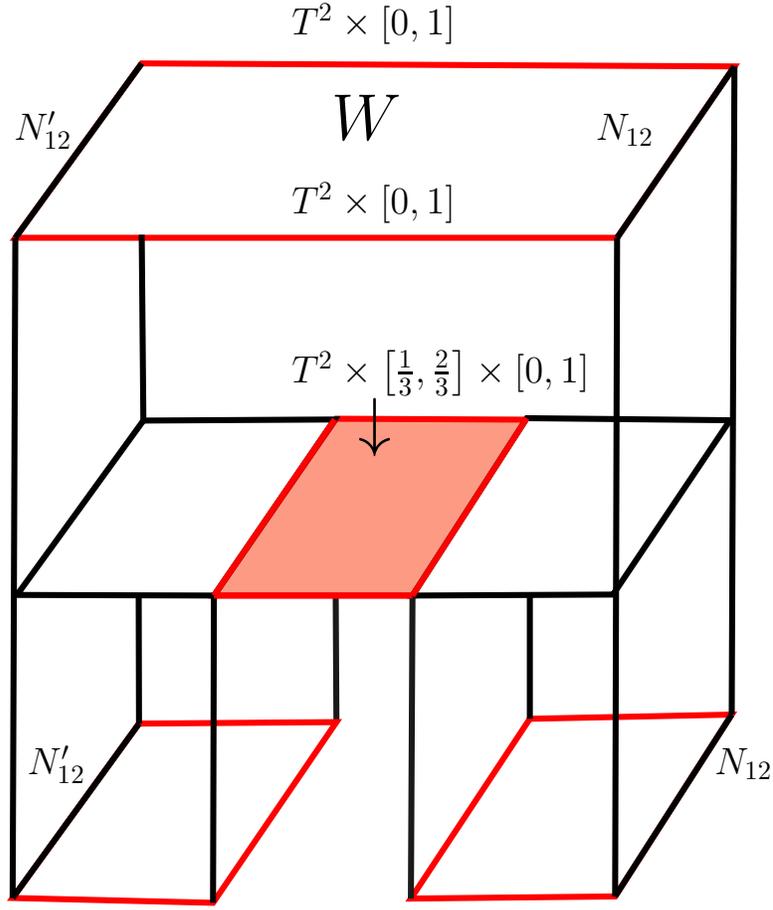}
	\put(-55,290){\large{$N_{12}$}}
	\put(-10,50){\large{$N_{12}$}}
	\put(-275,290){\large{$N'_{12}$}}
	\put(-270,50){\large{$N'_{12}$}}
	\put(-155,290){\huge{$W$}}
	\put(-170,262){\large{$T^{2}\times [0,1]$}}
	\put(-170,330){\large{$T^{2}\times [0,1]$}}
	\put(-170,198){\large{$T^{2}\times \left[\frac{1}{3},\frac{2}{3}\right] \times [0,1]$}}
	\put(-145,175){\Huge{$\downarrow$}}
	\caption{``Triangle equality'' and well-definedness of $\eta$.}
	\label{fig:surgequiv6}
\end{figure}

\begin{figure}[htb]
	\centering
		\includegraphics[scale = .4]{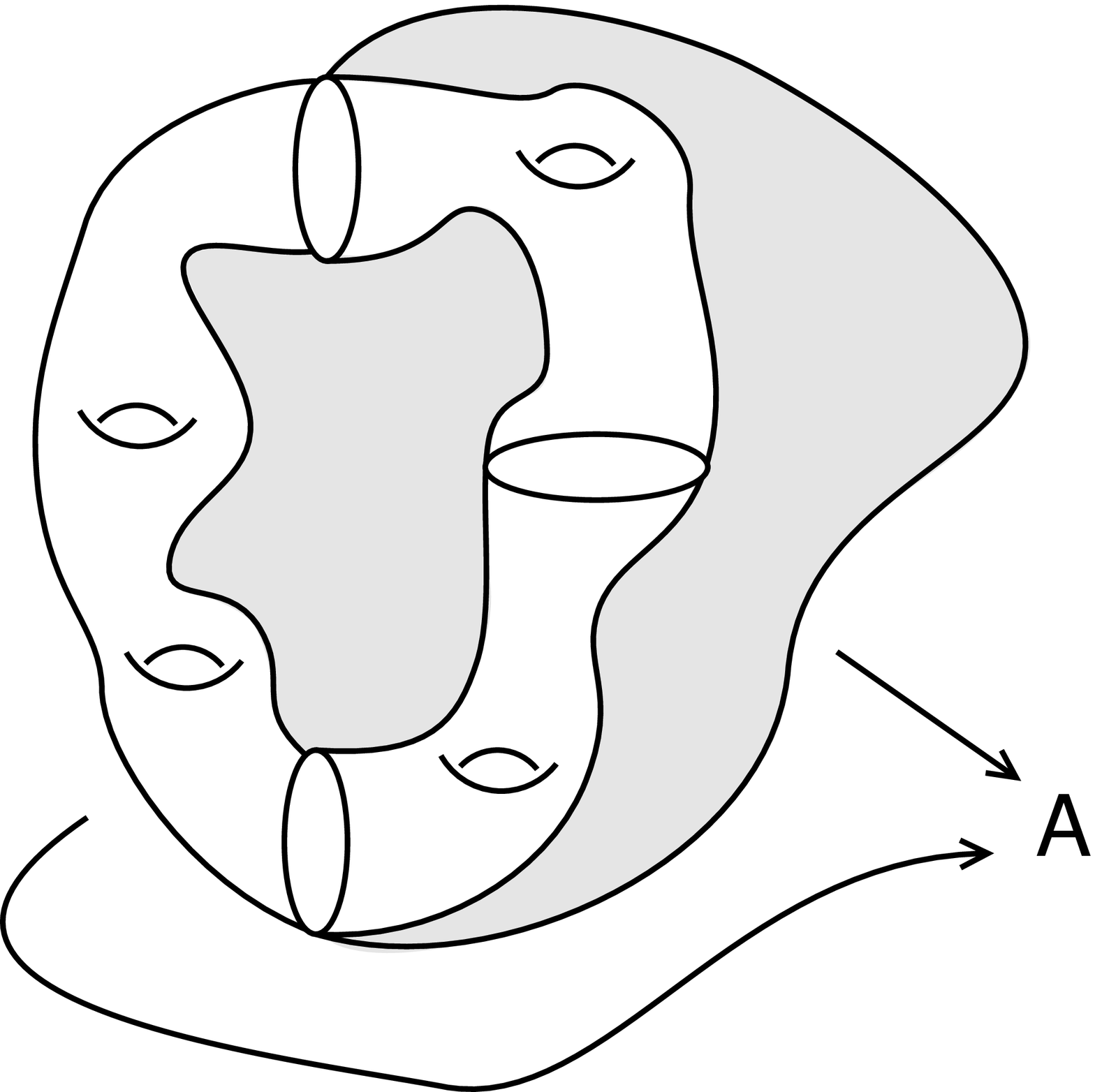}
	\put(-121,104){$N_{23}$}
	\put(-188,60){$N_{12}$}
	\put(-105,160){$N_{31}$}
	\caption{$N_{ij}$ glued together bord over $\Z_{2}$.}
	\label{fig:surgequiv5}
\end{figure} 


Let $W_{ij}$, and $N_{ij}$ denote the bordism and connecting $3$-manifolds between $S^{3}-K_{i}$ and $S^{3}-K_{j}$ for $1 \leq i,j \leq 3$ as above.  Also let $\overline{N_{ij}}$ and $\Phi_{ij}$ be as in the definition of $\eta(K_{i},K_{j})$.  

To prove well-definedness we must show that $\eta(K_{1},K_{2})$ is unchanged by any choice of connecting manifold.  Suppose there are two bordisms $(W_{12},\partial W_{12}, \Phi)$ and $(W'_{12},\partial W'_{12},\Phi')$ over $(D_{2p},\Z_{2})$ with connecting manifolds $N_{12}$ and $N'_{12}$.  Gluing $W_{12}$ together with $W'_{12}$ along their common boundaries $M_{1}=S^{3}-K_{1}$ and $M_{2}=S^{3}-K_{2}$ we see that $N_{12}\cup_{T^{2}\times \left\{0,1\right\}}N'_{12}$ bords over $\Z_{2}$.  Call this bordism $W$.  Up to bordism over $\Z_{2}$ we may assume that $\partial W= N_{12}\coprod N'_{12}\cup_{T^{2}\times \{0,1\}\times\{0,1\}}[(T^{2}\times [0,1])\times \{0,1\}]$ (see the top of Figure \ref{fig:surgequiv6}).  We may glue in a copy of $T^{2}\times \left[\frac{1}{3},\frac{2}{3}\right] \times [0,1]$ which shows that the disjoint union of $(\overline{N_{12}},\Phi_{12})$ and $\overline{N'_{21}},\Phi_{21})$ must also bord over $\Z_{2}$.  Of course Figure \ref{fig:surgequiv6} is just a rough diagram of this construction when thought of as a $5$-manifold.  That is $\Phi_{12,\ast}(\left[\overline{N_{12}}\right])+\Phi'_{21,\ast}(\left[\overline{N_{21}}\right])$.  Notice that $\overline{N_{21}}$ is just $\overline{N_{12}}$ with the reverse orientation but since we are working over $\Z_{2}$ the order does not matter.  That is, $$\Phi'_{21,\ast}(\left[\overline{N_{21}}\right])=-\Phi'_{12,\ast}(\left[\overline{N_{12}}\right])=\Phi'_{12,\ast}(\left[\overline{N_{12}}\right])\mod 2$$ and thus $\eta$ is invariant under the choice of bordism class $W_{12}$.

For the proof of the ``triangle equality'' we first obtain a bordism $W$ as in Figure \ref{fig:surgequiv5} by gluing all three $W_{ij}$'s along their common knot exterior boundaries.  In particular, the $3$-manifold obtained by gluing $N_{12}$, $N_{23}$, and $N_{31}$ together along their torus boundaries must bord over $\Z_{2}$.  But with a slight modification to the proof of well-definedness we obtain the relation $$\eta(K_{1},K_{2})+\eta(K_{2},K_{3})+\eta(K_{3},K_{1})=0 $$ in $\Omega_{3}(\Z_{2})\cong Z_{2}$ as desired. 

\end{proof}

So the bordism invariant $\omega(K,\rho)$ which is $\Zp$-valued may not be a complete invariant for surgery equivalence classes of $p$-colored knots.  However, if $(K_{1},\rho_{1})$ and $(K_{2},\rho_{2})$ are surgery equivalent based $p$-colored knots, then it is clear that $\omega(K_{1},\rho_{1})=\omega(K_{2},\rho_{2})$.  Recall that two $p$-colored knots are surgery equivalent if one may be obtained from the other by $\pm 1$-framed surgery on $S^{3}$ along an unlink $L=L_{1}\cup L_{2}$ with $[L_{i}]\in ker(\rho_{i})$ for $i=1,2$.  So the bordism over $(D_{2p},\Z_{2})$ is constructed by attaching $2$-handles along the components of $L_{1}$ and \textit{dual $2$-handles} along the components of $L_{2}$ to the $4$-manifold $(S^{3}-K_{1})\times [0,1]$.  Notice that the connecting manifold for this bordism is $T^{2}\times [0,1]$.  We have shown that surgery equivalent $p$-colored knots have the same bordism invariant.  The difficulty with the converse is indeed the connecting manifold.

We will now prove the main result.

\begin{proof}\textbf{\textit{Theorem \ref{thm:MainTheorem}}.}

By the discussion above, two surgery equivalent based $p$-colored knots have bordant exteriors over $(K(D_{2p}), K(\Z_{2}))$ where the $\Z_{2}=\left\langle t \right\rangle \subset D_{2p}=\left\langle s,t \; | \; s^{p}=t^{2}=stst=1 \right\rangle$.  If we assume that two \textit{$p$-colored knot exteriors} are bordant so that the connecting manifold is just the product space $T^{2}\times [0,1]$ then the converse is true.

Assume that $$(M_{1}=S^{3}-K_{1},f_{1})\sim_{(D_{2p},\Z_{2})} (M_{2}=S^{3}-K_{2},f_{2})$$ where the $f_{i}$ correspond to the coloring maps $\rho_{i}: \pi_{1}(S^{3}-K_{i})\rightarrow D_{2p}$ with bordism $(W^{4},\Omega)$.  Suppose further that $\partial W=(M_{1}\coprod -M_{2})\cup_{\partial T^{2}\times[0,1]}(T^{2}\times[0,1])$.  Take a smooth handle decomposition of $W$ relative to the boundary with no $0$ or $4$-handles and proceed in a similar way to the proof of Theorem 4.2 in \cite{CGO}.

We may ``trade'' $1$-handles for $2$-handles (see \cite[pages~6-7]{Kir} or \cite[Section 5.4]{GomSt}).  Since $(f_{1})_{\ast}:\pi_{1}(S^{3}-K_{1})\rightarrow D_{2p}$ is an epimorphism we may alter the attaching maps $c_{i}$ of the $2$-handles so that $\Phi_{\ast}(c_{i})=1$.  Thus the map $\Phi$ extends to the ``new'' $4$-manifold $W$ with no $1$-handles.  Since the $3$-handles may be thought of as upside down $1$-handles we may assume that $W$ is obtained from $(S^{3}-K_{1})\times [0,1]$ by attaching $2$-handles.  This implies that $M_{1}$ and $M_{2}$ are related by surgery along links in the kernel of $\rho_{i}$.  Now we must show that these links have $\pm 1$-framing and are unknotted.


Assume $\eta(K_{1},K_{2})=0$.  Then the connecting manifold is the product space $T^{2}\times [0,1]$.  So we may glue in a solid torus crossed with an interval to the boundary tori of $W$ and ``fill in'' the $M_{i}$ and the connecting manifold.  The result is a bordism between $S^{3}=(S^{3}-K_{1})\cup (S^{1}\times D^{2})$ and $S^{3}=(S^{3}-K_{2})\cup (S^{1}\times D^{2})$.  That is we have a surgery description of $S^{3}=(S^{3}-K_{2})\cup (S^{1}\times D^{2})$ consisting of a link $L$ in the complement of $K_{1}$.  We now appeal to \textit{Kirby's Theorem} to obtain the standard surgery description for $S^{3}$ by using only \textit{blow ups} and \textit{handle slides} and no \textit{blow downs}, consisting of a $\pm 1$-framed unlink.  Notice that by taking $K_{1}\subset S^{3}$ ``along for the ride'' when we do a handle slide we have only changed $K_{1}$ by an isotopy and so the resulting knot is surgery equivalent vacuously.  By a blow up we mean the addition of a single $\pm 1$-framed unknot away from the rest of the surgery diagram.  Since this unknot may be assumed to be in the kernel of $\rho_{1}$ and so this move is a surgery equivalence.  Hence we have shown that $(K_{1},\rho_{1})$ is surgery equivalent to $(K_{2},\rho_{2})$ if we assume that $\eta(K_{1},K_{2})=0$.

If $\eta(K_{1},K_{2})\neq 0$, then Proposition \ref{prop:EtaTriangle} implies that there are at most $2$ surgery classes of $p$-colored knots which have the same value of $\omega$.  As $\omega$ is $\Zp$-valued we have that there are no more than $2p$ possible equivalence classes.  Note that we have already seen that the connected sum of $k$ $(p,2)$-torus knots for $k=1,\ldots,p$ give a complete list of representatives for the $\Zp$-valued invariant $\omega$ and so the second statement of the proof follows from this.
\end{proof}

\end{document}